\documentclass[11pt]{article}
\usepackage{mathrsfs}
\usepackage{amsmath}
\usepackage{amsfonts}
\usepackage{enumerate}
\usepackage{setspace}
\usepackage{bbm}
\usepackage{amsthm}
\usepackage{amssymb}
\usepackage{xcolor}
\usepackage[title]{appendix}
\usepackage{cite}
\usepackage[colorlinks=true]{hyperref}
\usepackage[nameinlink]{cleveref}
\hypersetup{urlcolor=red, citecolor=blue}

\newtheorem{theorem}{Theorem}[section]
\newtheorem{proposition}[theorem]{Proposition}
\newtheorem{definition}[theorem]{Definition}
\newtheorem{lemma}[theorem]{Lemma}
\newtheorem{remark}[theorem]{Remark}

\crefname{subsection}{Subsection}{Subsections}
\crefname{theorem}{Theorem}{Theorems}
\crefname{example}{Example}{Examples}
\crefname{proposition}{Proposition}{Propositions}
\crefname{corollary}{Corollary}{Corollaries}
\crefname{definition}{Definition}{Definitions}
\crefname{lemma}{Lemma}{Lemmas}
\crefname{remark}{Remark}{Remarks}
\crefname{assumption}{Assumption}{Assumptions}
\crefname{problem}{Problem}{Problems}
\crefname{appendix}{Appendix}{Appendices}
\Crefname{appendix}{Appendix}{Appendices}

\hypersetup{
	hidelinks,
	colorlinks=true,
	linkcolor=red,
	citecolor=blue,
	urlcolor = black
}

\begin{document}

\title{Viscosity solutions of systems of second-order integral-partial differential equations with nonlinear jump terms without monotonicity conditions under infinite L\'{e}vy measures}
\author{
  Sa\"{i}d Hamad\`{e}ne$^{1}$,
  Hongyu Shi$^{2,}$\thanks{
    Corresponding author.
    E-mail:
    \href{mailto:202112026@mail.sdu.edu.cn}
         {\texttt{202112026@mail.sdu.edu.cn}}
  }\hphantom{$^{*}$},
  Zhen Wu$^{3}$
  \\[1.2ex]
  {\small $^{1}$Laboratoire Manceau de Math\'ematiques, Le Mans Universit\'e, Le Mans 72085, France}
  \\[0.8ex]
  {\small $^{2}$School of Mathematics, Shandong University, Jinan 250100, China}
  \\[0.8ex]
  {\small $^{3}$School of Mathematics, State Key Laboratory of Cryptography and Digital Economy Security,}
  \\
  {\small Shandong University, Jinan 250100, China}
}
\date{}  
\maketitle


\begin{abstract}
    This paper studies viscosity solutions of a system of second-order integral-partial differential equations with nonlinear jump terms, in which the jump increments are first transformed by a Lipschitz continuous function and then integrated against the L\'{e}vy measure. In contrast to the classical setting, no monotonicity conditions are imposed on the generator with respect to the jump component, and this function is allowed to have nonmonotone components and a nonzero value at the origin. We construct a viscosity solution by means of the connection between integral-partial differential equations and backward stochastic differential equations with Poisson jumps. When the L\'{e}vy measure is infinite, the jump component of the associated backward stochastic differential equation cannot be represented directly in terms of the constructed solution, and this difficulty is overcome by first truncating the L\'{e}vy measure to a finite measure and then passing to the limit. On this basis, we introduce an auxiliary integral-partial differential equation driven by a prescribed continuous field and thereby prove the existence and uniqueness of viscosity solutions to the original system. This result extends the known one for linear jump terms to the case of nonlinear jump terms. Specific types of these equations are the Hamilton--Jacobi--Bellman equations satisfied by value functions of firms which face Poisson shocks in an asymmetric manner.
\end{abstract}

\section{Introduction}
\label{sec1}
We consider the second-order integral-partial differential equations (IPDEs) of the following form:
\begin{equation}
    \label{eq001}
    \left\{
    \begin{aligned}
        &-\partial_tu_i(t,x)-\mathcal L u_i(t,x)- \mathcal Ku_i(t,x)\\
        &\qquad -h_i\left(t,x,u(t,x),\big(D_xu_i(t,x)\big)^\top \sigma(t,x),\mathcal B_i^\Xi u_i(t,x)\right)=0, \quad (t,x)\in[0,T]\times \mathbb{R}^k,\\
        & u_i(T,x)=g_i(x), \quad i = 1,\cdots,m,
    \end{aligned}
    \right.
\end{equation}
where
\begin{equation*}
    \left\{
        \begin{aligned}
            &\mathcal L u_i(t,x) = b(t,x)^\top D_xu_i(t,x)+\frac{1}{2}\operatorname{Tr}\left(\sigma(t,x)\sigma^\top(t,x)D^2_{xx}u_i(t,x)\right),\\
            &\mathcal K u_i(t,x) = \int_E\left[u_i(t,x+\beta(t,x,e))-u_i(t,x)-\beta(t,x,e)^\top D_xu_i(t,x)\right]\lambda(\mathrm{d}e),\\
            & \mathcal B_i^\Xi u_i(t,x) = \int_E\gamma_i(t,x,e)\Xi_i\left(u_i(t,x+\beta(t,x,e))-u_i(t,x)\right)\lambda(\mathrm{d}e).
        \end{aligned}
    \right.
\end{equation*}
Here, $\lambda$ is a L\'{e}vy measure on $E:=\mathbb R^\ell\setminus\{0\}$ satisfying $\int_E(1\wedge|e|^2)\lambda(\mathrm de)<\infty$, with $\lambda(E)=\infty$ allowed, and $\Xi=(\Xi_1,\cdots,\Xi_m)^\top$ is a Lipschitz continuous function. Associated with \eqref{eq001} is the following multidimensional backward stochastic differential equation (BSDE) with Poisson jumps:
\begin{equation}\label{bsdeintro}
    Y_{s}^{t,x}=g(X_T^{t,x})+\int_s^Tf\left(r,X_r^{t,x},Y_r^{t,x},Z_{r}^{t,x},U_{r}^{t,x}\right)\mathrm dr-\int_s^T Z_{r}^{t,x}\mathrm dW_r-\int_s^T\int_EU_{r}^{t,x}\tilde\mu(\mathrm dr,\mathrm de),\quad s\in [0,T],
\end{equation}
where $W$ is a Brownian motion, $\tilde\mu$ is a compensated Poisson random measure, and $X^{t,x}$ is the solution of the stochastic differential equation (SDE) with Poisson jumps starting from $(t,x)$, i.e., 
\begin{equation*}
    X_s^{t,x}=x+\int_t^s b(r,X_r^{t,x})\mathrm dr +\int_t^s\sigma(r,X_r^{t,x})\mathrm dW_r+\int_t^s\int_E\beta(r,X_{r-}^{t,x},e)\tilde\mu(\mathrm dr,\mathrm de),s\in(t,T];X_s^{t,x}=x,s\in[0,t].
\end{equation*}
The generator of \eqref{bsdeintro} is given by $f_i(t,x,y,z_i,\zeta_i)=h_i\big(t,x,y,z_i,\int_E\gamma_i(t,x,e)\Xi_i(\zeta_i(e))\lambda(\mathrm de)\big)$, $i=1,\cdots,m$. Since randomness enters the terminal condition and the generator only through $X^{t,x}$, $Y^{t,x}$ has a Markovian structure (see \Cref{pro23}(ii)). More precisely, there exists a continuous function $u$ of polynomial growth such that
\begin{equation}
    \label{eq001_02}
    Y_s^{t,x} = u(s,X_s^{t,x}),\quad \forall s\in [t,T],\quad \mathbb{P}\text{-a.s.}
\end{equation}

This class of IPDEs was introduced by Barles, Buckdahn and Pardoux~\cite{BBP1997}, who considered only the case $\Xi_i(r)=r$, $i=1,\cdots,m$. For the IPDE \eqref{eq001} considered here, the operator $\mathcal B_i^\Xi$ first applies the nonlinear transformation $\Xi_i$ to the jump increment $u_i(t,x+\beta(t,x,e))-u_i(t,x)$ and then integrates against $\gamma_i(t,x,e)\lambda(\mathrm de)$. It therefore allows different responses to jump shocks of different signs and magnitudes, while in \cite{BBP1997}, the structure is of compensation type. This feature is particularly useful when dealing with the determination of the solution of the IPDE of Hamilton--Jacobi--Bellman type satisfied by the value functions of firms which face Poisson shocks in an asymmetric manner. For example, credit valuation provides a natural motivation for such asymmetric responses. In the Merton model~\cite{M1974}, creditors receive $\min(V_T,B)=B-(B-V_T)^+$ at maturity, where $V_T$ denotes the value of the firm at debt maturity and $B$ is the promised repayment. Thus, creditors bear the entire shortfall when the firm value falls below $B$, but do not share in gains above $B$. This payoff asymmetry is an important motivation for the principle of accounting conservatism: recognizing gains requires a higher degree of verifiability than recognizing losses~\cite{B1997,W2003,OPK2023}. Accordingly, when assessing value increments caused by jump shocks, creditors assign greater weight to negative increments than to positive ones, leading to the choice
\begin{equation}\label{poidsauts}
    \Xi_i(r)=w_{1,i}r^+-w_{2,i}r^-,\qquad w_{2,i}>w_{1,i}\ge0,\quad i=1,\cdots,m.
\end{equation}
In this paper, $\Xi$ is only required to be Lipschitz continuous, with no monotonicity or differentiability assumption on its components. Moreover, $\Xi(0)$ captures the effect of each jump that is independent of the size of the increment, such as a fixed disposal cost incurred by each shock, in which case $\Xi(0)\neq0$.

The well-posedness of IPDEs can be studied by constructing viscosity solutions through BSDEs with Poisson jumps. In the case $\Xi_i(r)=r$, Barles, Buckdahn and Pardoux~\cite{BBP1997} assumed $\gamma_i\ge0$ and that $h_i(t,x,y,z_i,q)$ is nondecreasing in its last argument $q$ to ensure monotonicity of the equation with respect to the jump increment $u_i(t,x+\beta(t,x,e))-u_i(t,x)$. Under these two assumptions, referred to in the literature as monotonicity conditions, they proved that $u$ of \eqref{eq001_02} is the unique viscosity solution of the associated IPDE. Further results established under such conditions can be found in~\cite{BHL2011,LLW2025,DQS2015,DRZ2021,LLLW2025,LW2026}. In the present paper, however, no monotonicity conditions are imposed on \eqref{eq001}. Indeed, even if $\gamma_i\ge0$ and $h_i$ is nondecreasing in its last argument, a nonmonotone $\Xi_i$ can cause $h_i(\cdots,\mathcal B_i^\Xi u_i)$ to lose its monotonicity with respect to the jump increment. It is therefore necessary to study the case without monotonicity conditions. Hamad\`{e}ne and Morlais~\cite{HM2016} removed these conditions when $\lambda(E)<\infty$, and Hamad\`{e}ne~\cite{H2016} subsequently extended the result to $\lambda(E)=\infty$ by truncating the L\'{e}vy measure. This approach was later applied to systems with interconnected obstacles~\cite{HMN2023,HMN2026}, with~\cite{HMN2026} requiring $\int_E(1\wedge|e|)\lambda(\mathrm de)<\infty$. These models, however, involve only linear integrals of jump increments and cannot capture responses involving nonlinear transformations of those increments. Although nonlinear dependence of the generator on the jump component $U^{t,x}$ also arises in utility maximization and dynamic risk measures based on BSDEs with Poisson jumps~\cite{B2006,M2010,KPZ2015,QS2013,BLY2022,SST2024,BCCS2024,BHK2023}, a viscosity solution theory for the corresponding IPDEs has yet to be established.

The main objective of this paper is to establish the existence and uniqueness of viscosity solutions to the second-order IPDE system \eqref{eq001}, without imposing monotonicity conditions on $\gamma_i$, $h_i$, or $\Xi_i$, without requiring $\Xi$ to be differentiable or to satisfy $\Xi(0)=0$, and allowing $\lambda(E)=\infty$. Two points are crucial. {\bf The first point} is that the spatial difference quotients of the function $u$ in \eqref{eq001_02} are locally bounded (see \Cref{pro31}). Beyond the standard assumptions, this property requires only a compatibility condition on the composition of $h_i$ with $\int_E\gamma_i(t,x,e)\Xi_i(0)\lambda(\mathrm de)$. This condition can hold even when $\Xi(0)\neq0$ and $\int_E(1\wedge|e|)\lambda(\mathrm de)=\infty$ (see \Cref{re21}). {\bf The second point} is the characterization of the jump component $U^{t,x}$ in terms of $u$ and $X^{t,x}$ (see \Cref{th32}):
\begin{equation}
    \label{eq001_03}
    U^{t,x}_s(e) = u\big(s,X_{s-}^{t,x}+\beta(s,X_{s-}^{t,x},e)\big)-u\big(s,X_{s-}^{t,x}\big),\quad \mathrm{d}\mathbb{P} \otimes \mathrm{d}\lambda \otimes \mathrm{d}t\text{-a.e. on } \Omega \times E \times [t,T].
\end{equation}
Specifically, we truncate $\lambda$ near the origin to obtain a finite measure $\lambda^n$ and then let $n\to\infty$. Convergence of the truncated equations yields the Markovian representation of $U^{t,x}$. Unlike~\cite{H2016}, we allow $\Xi(0)\neq0$ and $\int_E(1\wedge|e|)\lambda(\mathrm de)=\infty$ to hold simultaneously. Accordingly, rather than truncating the entire jump integral $\int_E\gamma_i(t,x,e)\Xi_i(\zeta_i(e))\lambda(\mathrm de)$, we retain the full per-jump contribution $\int_E\gamma_i(t,x,e)\Xi_i(0)\lambda(\mathrm de)$ and truncate only the remaining part. This gives estimates for the truncated equations that are uniform in $n$ (see \Cref{re33}(i)).

The two points above establish the spatial regularity of $u$. We therefore introduce a new notion of viscosity solution (see \Cref{de41}), which allows us to retain $\mathcal B_i^\Xi u_i$ in the viscosity inequalities instead of replacing it with $\mathcal B_i^\Xi\varphi$ as in~\cite{BBP1997}, where $\varphi$ is a test function. Next, replacing the last argument of $h_i(t,x,y,z_i,q)$ by a prescribed continuous field $\vartheta_i$ gives the auxiliary IPDE \eqref{eq4_08}. For any viscosity solution $\bar u$ of \eqref{eq001}, choosing $\vartheta_i=\mathcal B_i^\Xi\bar u_i$ connects $\bar u$ with the function $u$ in \eqref{eq001_02}. Since the field $\vartheta$ is only continuous, a Markovian representation of the form \eqref{eq001_03} for the corresponding jump component cannot be obtained directly. To address this difficulty, we establish joint convergence of the functions $u^n$ associated with the truncated equations with respect to the truncation parameter $n$ and the initial state $x_n$: as $n\to\infty$ and $x_n\to x$, $u^n(t,x_n)\to u(t,x)$. This resolves the difficulty and ultimately yields the existence and uniqueness of viscosity solutions to the IPDE \eqref{eq001} (see \Cref{th43}). To the best of our knowledge, this result has not previously been proved in such a general framework.

The remainder of the paper is organized as follows. \Cref{sec2} introduces the notation, basic assumptions, and properties of SDEs and BSDEs with Poisson jumps. \Cref{sec3} establishes the Markovian representation of the jump component $U^{t,x}$ under an infinite L\'{e}vy measure. \Cref{sec4} proves the main result, namely, the existence and uniqueness of viscosity solutions to the IPDE \eqref{eq001}. \Cref{sec5} concludes the paper. The \hyperref[appA]{Appendix} recalls the classical definition of a viscosity solution. 

\section{Preliminaries and Basic Assumptions}
\label{sec2}

\subsection{Notation}
\label{subse21}

Let $T>0$, and let $(\Omega,\mathcal F,\mathbb F,\mathbb P)$ be a complete filtered probability space, where $\mathbb F=(\mathcal F_t)_{0\le t\le T}$. This space carries a $d$-dimensional standard Brownian motion $W=(W_t)_{0\le t\le T}$ and a Poisson random measure $\mu$ on $[0,T]\times E$, where $E=\mathbb R^\ell\setminus\{0\}$, $\ell\ge1$, and $\mathcal E$ is the Borel $\sigma$-algebra on $E$. Assume that $W$ and $\mu$ are independent and that $\mathbb F$ is the usual augmentation of their natural filtration, so that it is complete and right-continuous. The intensity measure of $\mu$ is $\lambda(\mathrm de)\mathrm dt$, where $\lambda$ is a $\sigma$-finite measure on $(E,\mathcal E)$ satisfying $\int_E(1\wedge|e|^2)\lambda(\mathrm de)<\infty$. Note that we allow for $\lambda(E)=\infty$. For any $A\in\mathcal E$ satisfying $\lambda(A)<\infty$ and every $t\in [0,T]$, the compensated Poisson random measure $\tilde\mu$ is defined by $\tilde\mu(\omega,[0,t]\times A):=\mu(\omega,[0,t]\times A)-t\lambda(A)$. It is well known that, for each $A$, $\tilde\mu(\omega, [0,t]\times A)$ is a martingale. Let $\mathcal P$ and $\mathcal G$ denote, respectively, the predictable and progressively measurable $\sigma$-algebras on $\Omega\times[0,T]$ with respect to $\mathbb F$.

We denote the $k$-dimensional (resp.~$k\times m$-dimensional) Euclidean space by $\mathbb R^k$ (resp.~$\mathbb R^{k\times m}$), with norm $|\cdot|$ and inner product $\langle\cdot,\cdot\rangle$. The transpose of a matrix $D$ is denoted by $D^\top$. For a c\`adl\`ag process $\theta$, that is, a process with right-continuous paths having left limits, we write $\theta_{t-}:=\lim_{s\uparrow t}\theta_s$ for $t\in(0,T]$, and set $\theta_{0-}:=\theta_0$. For a given Hilbert space $\mathbb H$, a filtration $\mathbb F$, positive integers $a,b$, and $t\in[0,T)$, $C^{a,b}([t,T]\times\mathbb R^k;\mathbb H)$ denotes the space of all $\mathbb H$-valued functions $\phi$ with continuous derivatives up to order $a$ in time and continuous partial derivatives up to order $b$ in space. $L^2(E;\mathbb H)$ denotes the space of all $\mathbb H$-valued Borel measurable functions $\zeta$ satisfying $\|\zeta\|^2_{L^2(E;\mathbb H)}=\int_E|\zeta(e)|^2\lambda(\mathrm de)<\infty$. $\mathcal S^2_{\mathbb F}([t,T];\mathbb H)$ denotes the space of all $\mathbb H$-valued, $\mathcal G$-measurable c\`adl\`ag processes $Y$ satisfying $\mathbb E[\sup_{t\le s\le T}|Y_s|^2]<\infty$. $\mathcal H^2_{\mathbb F}([t,T];\mathbb H)$ denotes the space of all $\mathbb H$-valued, $\mathcal G$-measurable processes $Z$ satisfying $\mathbb E\int_t^T|Z_s|^2\mathrm ds<\infty$. $\mathcal M^2_{\mathbb F}([t,T];\mathbb H)$ denotes the space of all $\mathbb H$-valued, $\mathcal P\otimes\mathcal E$-measurable processes $U$ satisfying $\mathbb E\int_t^T\|U_s\|_{L^2(E;\mathbb H)}^2\mathrm ds<\infty$.

Furthermore, $\Pi(\mathbb H)$ denotes the class of deterministic functions $\varphi:[0,T]\times\mathbb R^k\to\mathbb H$ of polynomial growth, i.e., there exist constants $C,p\ge0$ such that $|\varphi(t,x)|\le C(1+|x|^p)$ for all $(t,x)\in[0,T]\times\mathbb R^k$. Its subclass of continuous functions is denoted by $\Pi^c(\mathbb H)$. Let $\mathcal U(\mathbb H)$ denote the subclass of $\Pi^c(\mathbb H)$ consisting of functions $\varphi$ for which there exist constants $C,p\ge0$ such that 
\begin{equation*}
  |\varphi(t,x)-\varphi(t,x')|\le C(1+|x|^p+|x'|^p)|x-x'| \text{ for all } (t,x,x')\in[0,T]\times \mathbb R^{2k}.
\end{equation*}

\subsection{SDE with Poisson Jumps}
\label{subse22}

For any $(t,x)\in[0,T]\times\mathbb R^k$, consider the following SDE with Poisson jumps:
\begin{equation}
    \label{eq002}
    \left\{
        \begin{aligned}
            &X_s^{t,x}=x+\int_t^s b(r,X_r^{t,x})\mathrm dr +\int_t^s\sigma(r,X_r^{t,x})\mathrm dW_r+\int_t^s\int_E\beta(r,X_{r-}^{t,x},e)\tilde\mu(\mathrm dr,\mathrm de),\quad s\in(t,T],\\
            &X_s^{t,x}=x,\qquad s\in[0,t],
        \end{aligned}
    \right.
\end{equation}
where $(b,\sigma):[0,T]\times\mathbb R^k\to \mathbb R^k\times \mathbb R^{k\times d}$ and $\beta:[0,T]\times\mathbb R^k\times E\to \mathbb R^k$ are functions satisfying the following assumptions, which ensure the well-posedness of \eqref{eq002}.
\medskip

\noindent \textbf{Assumption 1.}
\begin{enumerate}[(i)]
\item The functions $b$ and $\sigma$ are continuous on $[0,T]\times\mathbb R^k$. Moreover, there exists a constant $C_{b,\sigma}>0$ such that, for all $(t,x,x')\in[0,T]\times \mathbb R^{2k}$,
\begin{equation*}
    |b(t,x)-b(t,x')|+|\sigma(t,x)-\sigma(t,x')|
    \le C_{b,\sigma}|x-x'|.
\end{equation*}
\item For any $e\in E$, $\beta(\cdot,\cdot,e)$ is continuous on $[0,T]\times\mathbb R^k$ and there exists a constant $C_\beta>0$ such that, for all $(t,x,x',e)\in[0,T] \times \mathbb R^{2k}\times E$,
\begin{equation*}
    |\beta(t,x,e)|\le C_\beta(1\wedge|e|)\mbox{ and } |\beta(t,x,e)-\beta(t,x',e)|\le C_\beta(1\wedge|e|)|x-x'|.
\end{equation*}
\end{enumerate}

By Theorem 2.1 and Lemma 2.1 of~\cite{FK1985}, the SDE \eqref{eq002} admits, under Assumption 1, a unique solution $X^{t,x} \in \mathcal S^2_{\mathbb F}([0,T];\mathbb R^k)$. Moreover, for any $p\ge2$, there exists a constant $C_p>0$ such that, for all $0\le t\le s\le T$ and $x,x'\in\mathbb R^k$,
\begin{equation}
    \label{eq003}
    \left\{
        \begin{aligned}
            &\mathbb E\left[\sup_{r\in[t,s]}|X_r^{t,x}-x|^p\right]\le C_p(s-t)(1+|x|^p),\\
            & \mathbb E\left[\sup_{r\in[t,s]}\left|X_r^{t,x}-X_r^{t,x'}-(x-x')\right|^p\right]\le C_p(s-t)|x-x'|^p. 
        \end{aligned}
    \right. 
\end{equation}

\subsection{Markovian BSDE with Poisson Jumps}
\label{subse23}
This subsection introduces the BSDE associated with IPDE \eqref{eq001}. The precise connection will be established in \Cref{sec4}. Here we present the basic formulation and properties of the BSDE. Let $h = (h_1,\cdots,h_m)^\top:[0,T]\times\mathbb R^k\times\mathbb R^m \times\mathbb R^{m\times d}\times\mathbb R^m \to\mathbb R^m$, $\gamma = (\gamma_1,\cdots,\gamma_m)^\top:[0,T]\times\mathbb R^k\times E \to \mathbb R^m$, and $\Xi=(\Xi_1,\cdots,\Xi_m)^\top:\mathbb R\to\mathbb R^m$ be measurable functions. The component functions $h_i:[0,T]\times\mathbb R^k\times\mathbb R^m \times\mathbb R^{1\times d}\times\mathbb R \to\mathbb R, i=1,\cdots,m$, are related to $h$ by
\begin{equation*}
    \begin{aligned}
        &h(t,x,y,z,q) = \big(h_1(t,x,y,z_1,q_1),\cdots,h_m(t,x,y,z_m,q_m)\big)^\top,\\
        &\qquad\qquad\qquad\qquad\qquad\qquad\qquad\qquad\qquad\qquad \forall (t,x,y,z,q) \in [0,T]\times\mathbb R^k\times\mathbb R^m \times\mathbb R^{m\times d}\times \mathbb R^m,
    \end{aligned}
\end{equation*}
where $z_i$ is the $i$-th row of the matrix $z$, and $q_i$ is the $i$-th component of the vector $q$. For any $(t,x)\in[0,T]\times\mathbb R^k$, let $X^{t,x}$ be the solution of \eqref{eq002}, and consider the following BSDE with Poisson jumps:
\begin{equation}
    \label{eq004}
    Y_{s}^{t,x}=g(X_T^{t,x})+\int_s^Tf\left(r,X_r^{t,x},Y_r^{t,x},Z_{r}^{t,x},U_{r}^{t,x}\right)\mathrm dr-\int_s^T Z_{r}^{t,x}\mathrm dW_r-\int_s^T\int_EU_{r}^{t,x}\tilde\mu(\mathrm dr,\mathrm de),\quad s\in [0,T],
\end{equation}
where $g = (g_1,\cdots,g_m)^\top:\mathbb R^k\to \mathbb R^m$ and $f= (f_1,\cdots,f_m)^\top:[0,T]\times\mathbb R^k\times\mathbb R^m \times\mathbb R^{m\times d}\times L^2(E;\mathbb R^m) \to\mathbb R^m$ are measurable functions. The function $f$ and its component functions $f_i:[0,T]\times\mathbb R^k\times\mathbb R^m \times\mathbb R^{1\times d}\times L^2(E;\mathbb R) \to\mathbb R$, $i=1,\cdots,m$, are given by
\begin{equation}
    \label{eq004_1}
    \begin{aligned}
        &f(t,x,y,z,\zeta)=\big(f_1(t,x,y,z_1,\zeta_1),\cdots,f_m(t,x,y,z_m,\zeta_m)\big)^\top\\
        &\!=\!\left(h_1\!\Big(t,x,y,z_1,\int_E\!\gamma_1(t,x,e)\Xi_1(\zeta_1(e))\lambda(\mathrm de)\Big),\cdots,h_m\!\Big(t,x,y,z_m,\int_E\!\gamma_m(t,x,e)\Xi_m(\zeta_m(e))\lambda(\mathrm de)\Big)\right)^\top\!,\\
        &\qquad\qquad\qquad\qquad\qquad\qquad\qquad\qquad\qquad \forall (t,x,y,z,\zeta) \in [0,T]\times\mathbb R^k\times\mathbb R^m \times\mathbb R^{m\times d}\times L^2(E;\mathbb R^m),
    \end{aligned}
\end{equation}
where $z_i$ is the $i$-th row of the matrix $z$, and $\zeta_i$ is the $i$-th component of $\zeta$. The following assumptions are imposed.
\medskip

\noindent \textbf{Assumption 2.}
\begin{enumerate}[(i)]
\item For any $(y,z,q)\in\mathbb R^m\times\mathbb R^{m\times d}\times\mathbb R^m$, the function $h(\cdot,\cdot,y,z,q)$ is continuous on $[0,T]\times\mathbb R^k$. There exist constants $C_h,p_h>0$ such that, for all $(t,x,x',y,y',z,z',q,q')\in[0,T]\times\mathbb R^{2k}\times\mathbb R^{2m} \times\mathbb R^{m\times d}\times\mathbb R^{m\times d}\times\mathbb R^{2m}$,
\begin{equation*}
    \left\{
        \begin{aligned}
            &|h(t,x,y,z,q)-h(t,x,y',z',q')|\le C_h\bigl(|y-y'|+|z-z'|+|q-q'|\bigr),\\[3pt]
            &|h(t,x,y,z,q)-h(t,x',y,z,q)|\le C_h(1+|x|^{p_h}+|x'|^{p_h})|x-x'|.
        \end{aligned}
    \right.
\end{equation*}
\item For any $e\in E$, the function $\gamma(\cdot,\cdot,e)$ is continuous on $[0,T]\times\mathbb R^k$. There exist constants $C_\gamma,p_\gamma>0$ such that, for all $(t,x,x',e)\in[0,T]\times\mathbb R^{2k}\times E$,
\begin{equation*}
    |\gamma(t,x,e)|\le C_\gamma(1\wedge|e|),\quad |\gamma(t,x,e)-\gamma(t,x',e)|\le C_\gamma(1\wedge|e|)(1+|x|^{p_\gamma}+|x'|^{p_\gamma})|x-x'|.
\end{equation*}
\item The function $\Xi$ is Lipschitz, i.e., there exists a constant $C_\Xi>0$ such that, for all $r,r'\in\mathbb R$,
\begin{equation*}
    |\Xi(r)-\Xi(r')|\le C_\Xi|r-r'|.
\end{equation*}
Moreover, for any $i=1,\cdots,m$ and $(t,x)\in[0,T]\times\mathbb R^k$, $\int_E|\gamma_i(t,x,e)\Xi_i(0)|\lambda(\mathrm de)<\infty$.
\item The function $g$ is continuous on $\mathbb R^k$ and there exist constants $C_g,p_g>0$ such that, for all $x,x'\in\mathbb R^{k}$,
\begin{equation*}
    |g(x)|\le C_g(1+|x|^{p_g}),\quad|g(x)-g(x')|\le C_g(1+|x|^{p_g}+|x'|^{p_g})|x-x'|.
\end{equation*}
\end{enumerate}

We further require the following assumptions on $h_i$, $\gamma_i$, and $\Xi_i$, $i=1,\cdots,m$:
\medskip

\noindent \textbf{Assumption 3.}
\begin{enumerate}[(i)]
\item For any $i = 1,\cdots,m$, $\int_E\gamma_i(\cdot,\cdot,e)\Xi_i(0)\lambda(\mathrm de) \in \Pi^c(\mathbb R)$.
\item There exist constants $\tilde C,\tilde p>0$ such that, for all $(t,x,x',y,z,q)\in[0,T]\times\mathbb R^{2k}\times\mathbb R^{m} \times\mathbb R^{1\times d}\times\mathbb R$,
\begin{equation*}
    \begin{aligned}
        &\quad \text{ }\left|h_i\left(t,x,y,z,q+\int_E\gamma_i(t,x,e)\Xi_i(0)\lambda(\mathrm de)\right)-h_i\left(t,x',y,z,q+\int_E\gamma_i(t,x',e)\Xi_i(0)\lambda(\mathrm de)\right)\right|\\
        &\le \tilde C(1+|x|^{\tilde p}+|x'|^{\tilde p})|x-x'|, \quad i=1,\cdots,m. 
    \end{aligned}
\end{equation*}
\end{enumerate}

\begin{remark}
    \label{re21}
    Under Assumption 2, either $\Xi(0) = 0$ or $\int_E(1\wedge|e|)\lambda(\mathrm de)<\infty$ is sufficient for Assumption 3 to hold automatically, but the converse is false. For example, take $m=k=\ell=1$ and
    \begin{equation*}
        \lambda(\mathrm de)=\mathbf 1_{(0,1)}(e)\frac{\mathrm de}{e^2},\quad \Xi(r)=1+r,\quad \gamma(t,x,e)=1\wedge|e|^2,\quad h(t,x,y,z,q)=q.
    \end{equation*}
    It is straightforward to verify that $\int_E(1\wedge|e|^2)\lambda(\mathrm de)=1<\infty$ and that Assumptions 2 and 3 hold, whereas $\Xi(0) = 1$ and $\int_E(1\wedge|e|)\lambda(\mathrm de)=\infty$.
\end{remark}

\begin{remark}
    \label{re22}
    Assumptions 2 and 3 imply the following properties:
    \begin{enumerate}[(i)]
    \item For any $i = 1,\cdots,m$ and $\zeta\in L^2(E;\mathbb R)$, $\int_E\gamma_i(\cdot,\cdot,e)\Xi_i(\zeta(e))\lambda(\mathrm de) \in \Pi^c(\mathbb R)$.
    \item $f$ is continuous on $[0,T]\times\mathbb R^k\times\mathbb R^m \times\mathbb R^{m\times d}\times L^2(E;\mathbb R^m)$, and $f(\cdot,\cdot,0,0,{\bf 0})\in \Pi^c(\mathbb R^m)$.
    \item There exists a constant $C_f>0$ such that, for all $(t,x,y,y',z,z',\zeta,\zeta')\in[0,T]\times\mathbb R^{k}\times\mathbb R^{2m} \times\mathbb R^{m\times d}\times\mathbb R^{m\times d}\times L^2(E;\mathbb R^m)\times L^2(E;\mathbb R^m)$,
    \begin{equation*}
        |f(t,x,y,z,\zeta)-f(t,x,y',z',\zeta')|\le C_f\bigl(|y-y'|+|z-z'|+\|\zeta-\zeta'\|_{L^2(E;\mathbb R^m)}\bigr).
    \end{equation*}
    \end{enumerate}
\end{remark}

Under Assumptions 1--3, Theorem 2.1 of~\cite{BBP1997} guarantees that \eqref{eq004} admits a unique solution $(Y^{t,x},Z^{t,x},U^{t,x}) \in \mathcal S^2_{\mathbb F}([0,T];\mathbb R^m)\times \mathcal H^2_{\mathbb F}([0,T];\mathbb R^{m\times d}) \times \mathcal M^2_{\mathbb F}([0,T];\mathbb R^m)$. For the reader's convenience, we recall the known Markovian properties of $(Y^{t,x},Z^{t,x},U^{t,x})$ below. Detailed proofs can be found in the references cited in \Cref{re24} and \Cref{re26}. 

\begin{proposition}
    \label{pro23}
    Under Assumptions 1--3, the solution $(Y^{t,x},Z^{t,x},U^{t,x})$ of \eqref{eq004} has the following properties:
    \begin{enumerate}[(i)]
    \item For any $p\geq 2$, there exists a positive constant $K(p,T,C_f)$ depending only on $p$, $T$, and the constant $C_f$ in \Cref{re22}(iii) such that
    \begin{equation*}
        \begin{aligned}
            &\quad \text{ } \mathbb E\left[\sup_{0\le s\le T}|Y_s^{t,x}|^p+\left(\int_0^T|Z_s^{t,x}|^2\mathrm ds\right)^{\frac{p}{2}}+\left(\int_0^T\|U_s^{t,x}\|_{L^2(E;\mathbb R^m)}^2\mathrm ds\right)^{\frac{p}{2}}\right]\\
            &\leq K(p,T,C_f) \mathbb E\left[\big|g(X_T^{t,x})\big|^p+\left(\int^T_0\big|f(s,X_s^{t,x},0,0,{\bf 0})\big|\mathrm ds\right)^p\right].
        \end{aligned}
    \end{equation*}
    \item There exists a function $u=(u_1,\cdots,u_m)^\top\in \Pi^c(\mathbb R^m)$ such that, for any $(t,x)\in [0,T]\times\mathbb R^k$,
    \begin{equation*}
        Y_s^{t,x} = u(s,X_s^{t,x}),\quad \forall s\in [t,T],\quad \mathbb{P}\text{-a.s.}, 
    \end{equation*}
    and $u(T,x) = g(x)$.
    \end{enumerate}
\end{proposition}

\begin{remark}
    \label{re24}
    Some comments on this proposition are in order.
    \begin{enumerate}[(i)]
    \item The conclusions of \Cref{pro23} remain valid for any generator $f$ satisfying only the properties in \Cref{re22}(ii) and~(iii), regardless of the specific integral representation in \eqref{eq004_1}.
    \item Theorem 2 of~\cite{EFO2017} ensures that, for any $p\geq 2$, $\mathbb E[\sup_{0\le s\le T}|Y_s^{t,x}|^p]<\infty$. Combining this with Lemma 4 of~\cite{EFO2017} directly yields \Cref{pro23}(i).
    \item The conclusion of \Cref{pro23}(ii) follows from Proposition 2.5 of~\cite{BBP1997} or Proposition 2.2 of~\cite{H2016}. Although \eqref{eq004_1} introduces the nonlinear function $\Xi$, the proof relies only on the properties in \Cref{re22}(ii) and~(iii), so \Cref{pro23}(ii) remains valid.
    \end{enumerate}
\end{remark}

\begin{proposition}
    \label{pro25}
    Under Assumptions 1--3, if $\lambda(E)<\infty$, then, for any $(t,x)\in [0,T]\times\mathbb R^k$,
    \begin{equation*}
        U^{t,x}_s(e) = u\big(s,X_{s-}^{t,x}+\beta(s,X_{s-}^{t,x},e)\big)-u\big(s,X_{s-}^{t,x}\big),\quad \mathrm{d}\mathbb{P} \otimes \mathrm{d}\lambda \otimes \mathrm{d}t\text{-a.e. on } \Omega \times E \times [t,T], 
    \end{equation*}
    where the deterministic function $u$ is given by \Cref{pro23}(ii).
\end{proposition}

\begin{remark}
    \label{re26}
    Proposition 3.1 of~\cite{HM2016} establishes \Cref{pro25}. Its proof characterizes $U^{t,x}$ using the continuous Markovian representation in \Cref{pro23}(ii), which yields
    \begin{equation*}
        \begin{aligned}
            &\int_t^s\int_EU_{r}^{t,x}\mu(\mathrm dr,\mathrm de) = \sum_{t<r\le s}\bigl(Y_r^{t,x}-Y_{r-}^{t,x}\bigr) = \sum_{t<r\le s}\left[u(r,X_r^{t,x})-u(r,X_{r-}^{t,x})\right]\\
            &\qquad\qquad\qquad\qquad\text{ }\text{ } = \int_t^s\int_E \left[u\bigl(r,X_{r-}^{t,x}+\beta(r,X_{r-}^{t,x},e)\bigr)-u(r,X_{r-}^{t,x})\right]\mu(\mathrm dr,\mathrm de),
        \end{aligned}
    \end{equation*}
    as in (3.2) and (3.3) of~\cite{HM2016}. Assumptions 2 and 3 are used primarily to establish the Markovian representation of $Y^{t,x}$. Thus, the hypotheses of \Cref{pro25} can be weakened to require only Assumption 1, the well-posedness of \eqref{eq004}, and \Cref{pro23}(ii). The specific integral representation in \eqref{eq004_1} and the nonlinear function $\Xi$ do not affect the identities above and play no role in this proof. 
\end{remark}

\section{Markovian Representation of the Jump Component under an Infinite L\'{e}vy Measure}
\label{sec3}
This section establishes the Markovian representation of the jump component in \Cref{pro25} when $\lambda(E)=\infty$. Recall that $(Y^{t,x},Z^{t,x},U^{t,x})$ is the unique solution of \eqref{eq004} and that there exists a deterministic function $u=(u_1,\cdots,u_m)^\top\in \Pi^c(\mathbb R^m)$ such that, for any $(t,x)\in [0,T]\times\mathbb R^k$, $u(T,x) = g(x)$ and
\begin{equation}
    \label{eq006}
    Y_s^{t,x} = u(s,X_s^{t,x}),\quad \forall s\in [t,T],\quad \mathbb{P}\text{-a.s.} 
\end{equation}
We first study the spatial regularity of the function $u$ in \eqref{eq006}.
\begin{proposition}
    \label{pro31}
    Under Assumptions 1--3, $u \in \mathcal{U}(\mathbb R^m)$, i.e., $u$ is continuous and there exist constants $C,p\ge 0$ such that 
    \begin{equation}\label{croipolyu}|u(t,x)-u(t,x')|\le C(1+|x|^p+|x'|^p)|x-x'|\end{equation}for any $(t,x,x')\in [0,T]\times \mathbb R^{k+k}$.
\end{proposition}
\begin{proof}
    The generator $f$ defined in \eqref{eq004_1} satisfies a useful spatial estimate. By Assumption~2(i)--(iii), Assumption~3(ii), and the Cauchy--Schwarz inequality, for any $(t,x,x',y)\in[0,T]\times\mathbb R^{2k}\times\mathbb R^{m}$, $z = (z_1,\cdots,z_m)^\top \in \mathbb R^{m\times d}$, and $\zeta = (\zeta_1,\cdots,\zeta_m)^\top \in L^2(E;\mathbb R^m)$,
    \begin{equation*}
        \begin{aligned}
            &\quad \text{ }|f_i(t,x,y,z_i,\zeta_i)-f_i(t,x',y,z_i,\zeta_i)|\\
            &\le\left|h_i\left(t,x,y,z_i,\int_E\gamma_i(t,x,e)[\Xi_i(\zeta_i(e))-\Xi_i(0)]\lambda(\mathrm de)+\int_E\gamma_i(t,x,e)\Xi_i(0)\lambda(\mathrm de)\right)\right.\\
            &\qquad\left.-h_i\left(t,x',y,z_i,\int_E\gamma_i(t,x,e)[\Xi_i(\zeta_i(e))-\Xi_i(0)]\lambda(\mathrm de)+\int_E\gamma_i(t,x',e)\Xi_i(0)\lambda(\mathrm de)\right)\right|\\
            &\qquad+C_h\int_E|\gamma_i(t,x,e)-\gamma_i(t,x',e)||\Xi_i(\zeta_i(e))-\Xi_i(0)|\lambda(\mathrm de)\\
            &\le\tilde C(1+|x|^{\tilde p}+|x'|^{\tilde p})|x-x'|\\
            &\quad+C_hC_\gamma C_\Xi\left(\int_E(1\wedge|e|^2)\lambda(\mathrm de)\right)^{\frac{1}{2}} (1+|x|^{p_\gamma}+|x'|^{p_\gamma})|x-x'|\|\zeta_i\|_{L^2(E;\mathbb R)},\quad i = 1,\cdots,m.
        \end{aligned}
    \end{equation*}
    Consequently,
    \begin{equation}
        \label{eq007}
        |f(t,x,y,z,\zeta)-f(t,x',y,z,\zeta)|\leq C\big(1+\|\zeta\|_{L^2(E;\mathbb R^m)}\big)(1+|x|^p+|x'|^p)|x-x'|,
    \end{equation}
    where $C,p>0$ are constants. Throughout the proof, the constants $C$ and $p$ may vary from line to line.

    For $(t,s,x,x')\in [0,T]\times[0,T]\times\mathbb R^{2k}$, define
    \begin{equation*}
        \left\{
            \begin{aligned}
                &\Delta X_s:=X_s^{t,x}-X_s^{t,x'},\quad \Delta Y_s:=Y_s^{t,x}-Y_s^{t,x'},\quad \Delta Z_s:=Z_s^{t,x}-Z_s^{t,x'},\quad \Delta U_s:=U_s^{t,x}-U_s^{t,x'},\\
                & \Delta f_s:=f(s,X_s^{t,x},Y_s^{t,x},Z_s^{t,x},U_s^{t,x}) -f(s,X_s^{t,x'},Y_s^{t,x'},Z_s^{t,x'},U_s^{t,x'}),\quad \Delta g_T := g(X_T^{t,x})-g(X_T^{t,x'}).
            \end{aligned}
        \right.
    \end{equation*}
    Applying It\^o's formula to $|\Delta Y_s|^2$ yields
    \begin{equation}
        \label{eq008}
        \mathbb E|\Delta Y_s|^2+\mathbb E\int_s^T\Big(|\Delta Z_r|^2+\|\Delta U_r\|_{L^2(E;\mathbb R^m)}^2\Big)\mathrm dr=\mathbb E|\Delta g_T|^2 +2\mathbb E\int_s^T\langle\Delta Y_r,\Delta f_r\rangle\mathrm dr.
    \end{equation}
    Combining \Cref{re22}(iii), \eqref{eq007}, and Young's inequality gives
    \begin{equation}
        \label{eq009}
        \begin{aligned}
            2\langle\Delta Y_r,\Delta f_r\rangle&\le C|\Delta Y_r|^2 +\frac{1}{2}\left(|\Delta Z_r|^2+\|\Delta U_r\|_{L^2(E;\mathbb R^m)}^2\right)\\
            &\quad+C(1+|X_r^{t,x}|^p+|X_r^{t,x'}|^p)^2|\Delta X_r|^2 \left(1+\|U_r^{t,x'}\|_{L^2(E;\mathbb R^m)}^2\right).
        \end{aligned}
    \end{equation}
    Moreover, for any $q\geq2$, \eqref{eq003} implies
    \begin{equation}
        \label{eq010}
        \mathbb E\left[\sup_{0\le s\le T}(|X_s^{t,x}|^q+|X_s^{t,x'}|^q)\right]\le C_q(1+|x|^q+|x'|^q),\quad \mathbb E\left[\sup_{0\le s\le T}|\Delta X_s|^q\right]\le C_q|x-x'|^q.
    \end{equation}
    By Assumption~2(iv), \Cref{re22}(ii), \Cref{pro23}(i), and \eqref{eq010},
    \begin{equation}
        \label{eq011_1}
        \mathbb E\left[\left(\int_0^T\|U_s^{t,x'}\|_{L^2(E;\mathbb R^m)}^2\mathrm ds\right)^{\frac{q}{2}}\right]\leq C_q (1+|x'|^{\varrho_q}),
    \end{equation}
    where $C_q,\varrho_q>0$ are constants. Hence, Assumption~2(iv) and \eqref{eq008}--\eqref{eq011_1} imply
    \begin{equation}
        \label{eq011}
        \mathbb E|\Delta Y_s|^2+\frac{1}{2}\mathbb E\int_s^T\left(|\Delta Z_r|^2+\|\Delta U_r\|_{L^2(E;\mathbb R^m)}^2\right)\mathrm dr\le C(1+|x|^p+|x'|^p)|x-x'|^2 +C\int_s^T\mathbb E|\Delta Y_r|^2\mathrm dr.
    \end{equation}
    In particular, setting $s=t$ in \eqref{eq011} and applying \eqref{eq006} and Gronwall's inequality yields
    \begin{equation}
        \label{eq012_1}
        |u(t,x)-u(t,x')| = |\Delta Y_t| \le C(1+|x|^p+|x'|^p)|x-x'|.
    \end{equation}
    Since $u\in \Pi^c(\mathbb R^m)$, \eqref{eq012_1} shows that $u\in \mathcal{U}(\mathbb R^m)$.
\end{proof}

We now establish the Markovian representation of $U^{t,x}$ for a general L\'{e}vy measure, allowing $\lambda(E)=\infty$.
\begin{theorem}
    \label{th32}
    Under Assumptions 1--3, for any $(t,x)\in [0,T]\times\mathbb R^k$,
    \begin{equation*}
        U^{t,x}_s(e) = u\big(s,X_{s-}^{t,x}+\beta(s,X_{s-}^{t,x},e)\big)-u\big(s,X_{s-}^{t,x}\big),\quad \mathrm{d}\mathbb{P} \otimes \mathrm{d}\lambda \otimes \mathrm{d}t\text{-a.e. on } \Omega \times E \times [t,T].
    \end{equation*}
\end{theorem}
\begin{proof}
    {\bf \underline{Step 1}: Truncation of the L\'{e}vy measure.} For any $n\geq 1$, define
    \begin{equation}
        \label{eq014_1}
        \left\{
            \begin{aligned}
                &E^n:=\bigg\{e\in E : |e|\ge \frac{1}{n}\bigg\},\qquad \lambda^n(\mathrm de):=\mathbbm{1}_{E^n}(e)\lambda(\mathrm de),\\
                &\mu^n(\mathrm ds,\mathrm de):=\mathbbm{1}_{E^n}(e)\mu(\mathrm ds,\mathrm de), \qquad\tilde\mu^n(\mathrm ds,\mathrm de):= \mu^n(\mathrm ds,\mathrm de)-\lambda^n(\mathrm de)\mathrm ds.
            \end{aligned}
        \right.
    \end{equation}
    Consider the corresponding SDE with Poisson jumps:
    \begin{equation}
        \label{eq012}
        \left\{
            \begin{aligned}
                &X_s^{n,t,x}=x+\int_t^s b(r,X_r^{n,t,x})\mathrm dr +\int_t^s\sigma(r,X_r^{n,t,x})\mathrm dW_r+\int_t^s\int_E\beta(r,X_{r-}^{n,t,x},e)\tilde\mu^n(\mathrm dr,\mathrm de),\text{ } s\in(t,T],\\
                &X_s^{n,t,x}=x,\qquad s\in[0,t].
            \end{aligned}
        \right.
    \end{equation}
    The SDE \eqref{eq012} is formulated with respect to $\mathbb{F}^n$, the usual augmentation of the natural filtration generated by the Brownian motion $W$ and the independent random measure $\mu^n$. Since $\mathbb{F}^n\subseteq \mathbb{F}$, Theorem 2.1 of~\cite{FK1985} ensures that, for any $n\geq 1$, \eqref{eq012} has a unique solution $X^{n,t,x} \in \mathcal S^2_{\mathbb F^n}([0,T];\mathbb R^k)\subseteq \mathcal S^2_{\mathbb F}([0,T];\mathbb R^k)$. The same observation applies to \eqref{eq013} below. For any $n\geq 1$, $\int_E(1\wedge|e|^2)\lambda^n(\mathrm de)\le\int_E(1\wedge|e|^2)\lambda(\mathrm de)<\infty$. Combining this bound with Assumption~1(ii) and \eqref{eq003} gives, for any $q\geq 2$ and $x,x'\in \mathbb R^k$,
    \begin{equation}
        \label{eq013_1}
        \left\{
            \begin{aligned}
                &\sup_{n\ge1}\mathbb E\left[\sup_{0\le s\le T}|X_s^{n,t,x}|^q\right]\le C_q(1+|x|^q),\\
                &\sup_{n\ge1}\mathbb E\left[\sup_{0\le s\le T}|X_s^{n,t,x}-X_s^{n,t,x'}|^q\right]\le C_q|x-x'|^q,
            \end{aligned}
        \right.
    \end{equation}
    where $C_q>0$ is a constant. In what follows, the constant $C_q>0$ may vary from line to line. We next introduce the BSDE with Poisson jumps corresponding to \eqref{eq004}:
    \begin{equation}
        \label{eq013}
        \begin{aligned}
            &Y_{s}^{n,t,x}=g(X_T^{n,t,x})+\int_s^Tf^n\left(r,X_r^{n,t,x},Y_r^{n,t,x},Z_{r}^{n,t,x},U_{r}^{n,t,x}\right)\mathrm dr\\
            &\qquad\qquad-\int_s^T Z_{r}^{n,t,x}\mathrm dW_r-\int_s^T\int_EU_{r}^{n,t,x}\tilde\mu^n(\mathrm dr,\mathrm de),\qquad s\in [0,T],
        \end{aligned}
    \end{equation}
    where, for any $(t,x,y)\in[0,T]\times\mathbb R^{k}\times\mathbb R^{m}$, $z = (z_1,\cdots,z_m)^\top \in \mathbb R^{m\times d}$, and $\zeta = (\zeta_1,\cdots,\zeta_m)^\top \in L^2(E;\mathbb R^m)$,
    \begin{equation}
        \label{eq019_2}
        \left\{
            \begin{aligned}
                &f^n(t,x,y,z,\zeta)=\big(f^n_1(t,x,y,z_1,\zeta_1),\cdots,f^n_m(t,x,y,z_m,\zeta_m)\big)^\top,\\
                & f_i^n(t,x,y,z_i,\zeta_i)\!=\!h_i\left(t,x,y,z_i,\int_{E}\gamma_i(t,x,e)[\Xi_i(\zeta_i(e))-\Xi_i(0)]\lambda^n(\mathrm de)+\int_E\gamma_i(t,x,e)\Xi_i(0)\lambda(\mathrm de)\right).
            \end{aligned}
        \right.
    \end{equation}
    As in \Cref{re22}(iii) and \eqref{eq007},
    \begin{equation}
        \label{eq014}
        \left\{
            \begin{aligned}
                &|f^n(t,x,y,z,\zeta)-f^n(t,x,y',z',\zeta')|\le C\bigl(|y-y'|+|z-z'|+\|\zeta-\zeta'\|_{L^2(E;\mathbb R^m)}\bigr),\\
                &|f^n(t,x,y,z,\zeta)-f^n(t,x',y,z,\zeta)|\leq C\big(1+\|\zeta\|_{L^2(E;\mathbb R^m)}\big)(1+|x|^p+|x'|^p)|x-x'|,\quad n\geq 1,
            \end{aligned}
        \right.
    \end{equation}
    where the constants $C,p>0$ are independent of $n$. Throughout the remainder of the proof, the constants $C$ and $p$ may vary from line to line and remain independent of $n$. For any $n\geq 1$, Theorem 2.1 of~\cite{BBP1997} ensures that \eqref{eq013} has a unique adapted solution $(Y^{n,t,x},Z^{n,t,x},U^{n,t,x})$. By \Cref{pro23}(ii), \Cref{re24}(i), and \eqref{eq014}, for each $n\geq 1$ there exists a deterministic function $u^n=(u^n_1,\cdots,u^n_m)^\top\in \Pi^c(\mathbb R^m)$ such that, for any $(t,x)\in [0,T]\times\mathbb R^k$, $u^n(T,x) = g(x)$ and
    \begin{equation}
        \label{eq016}
        Y_s^{n,t,x} = u^n(s,X_s^{n,t,x}),\quad \forall s\in [t,T],\quad \mathbb{P}\text{-a.s.} 
    \end{equation}
    Since the constants in \eqref{eq014} are independent of $n$, Assumption~2(iv), \Cref{re22}(ii), \Cref{pro23}(i), and \eqref{eq013_1} imply that, for any $q\geq 2$,
    \begin{equation}
        \label{eq019_1}
        \sup_{n\geq 1}\mathbb E\left[\left(\int_0^T\|U_s^{n,t,x}\|_{L^2(E;\mathbb R^m)}^2\mathrm ds\right)^{\frac{q}{2}}\right]\leq C_q (1+|x|^{\varrho_q}),
    \end{equation}
    where $\varrho_q>0$ is a constant. With \eqref{eq013_1}, \eqref{eq014}, and \eqref{eq019_1}, the argument leading to \eqref{eq008}--\eqref{eq012_1} in the proof of \Cref{pro31} yields
    \begin{equation}
        \label{eq018}
        \sup_{n\geq 1}|u^n(t,x)-u^n(t,x')| \le C(1+|x|^p+|x'|^p)|x-x'|,\quad (t,x,x')\in [0,T]\times\mathbb R^{2k}.
    \end{equation}
    Furthermore, since $\lambda^n(E)\le n^2\int_E(1\wedge|e|^2)\lambda(\mathrm de)<\infty$, \Cref{pro25} gives
    \begin{equation}
        \label{eq019}
        U^{n,t,x}_s(e) = u^n\big(s,X_{s-}^{n,t,x}+\beta(s,X_{s-}^{n,t,x},e)\big)-u^n\big(s,X_{s-}^{n,t,x}\big),\quad \mathrm{d}\mathbb{P} \otimes \mathrm{d}\lambda \otimes \mathrm{d}t\text{-a.e. on } \Omega \times E \times [t,T].
    \end{equation}

    \noindent {\bf \underline{Step 2}: Convergence of \eqref{eq012} and \eqref{eq013}.} For $n\geq 1$ and $(t,s,x,x')\in [0,T]\times[0,T]\times\mathbb R^{2k}$, define
    \begin{equation*}
        \Delta X^n_s:=X_s^{n,t,x}-X_s^{t,x},\quad \Delta Y^n_s:=Y_s^{n,t,x}-Y_s^{t,x},\quad \Delta Z^n_s:=Z_s^{n,t,x}-Z_s^{t,x},\quad \Delta U^n_s:=U_s^{n,t,x}-U_s^{t,x}.
    \end{equation*}
    We first consider convergence of the SDE with Poisson jumps. By \eqref{eq002} and \eqref{eq012},
    \begin{equation*}
        \begin{aligned}
            & \Delta X^n_s = \int_t^s\big[b(r,X_r^{n,t,x})-b(r,X_r^{t,x})\big]\mathrm dr +\int_t^s\big[\sigma(r,X_r^{n,t,x})-\sigma(r,X_r^{t,x})\big]\mathrm dW_r\\
            &\qquad\qquad+\int_t^s\int_{E^n}\big[\beta(r,X_{r-}^{n,t,x},e)-\beta(r,X_{r-}^{t,x},e)\big]\tilde\mu(\mathrm dr,\mathrm de)+\int_t^s\int_{E\setminus E^n}\beta(r,X_{r-}^{t,x},e)\tilde\mu(\mathrm dr,\mathrm de).
        \end{aligned}
    \end{equation*}
    For any $q\geq 2$, Assumption 1, the Burkholder--Davis--Gundy inequality, and Kunita's first inequality (see Proposition 6.1 of~\cite{FS2018} and Lemma 2.1.5 of~\cite{JP2011}) yield
    \begin{equation*}
        \begin{aligned}
            &\mathbb{E} \left[\sup_{0\leq s\leq T}|\Delta X^n_s|^q\right]\leq C_q\mathbb{E}\int_0^T|\Delta X^n_s|^q\mathrm ds+ C_q\mathbb{E}\left(\int^T_0\int_{E^n}|\beta(s,X_{s-}^{n,t,x},e)-\beta(s,X_{s-}^{t,x},e)|^2\lambda(\mathrm de)\mathrm ds\right)^{\frac{q}{2}}\\
            &\qquad\qquad\qquad\qquad\text{ } + C_q\mathbb{E}\int_0^T\int_{E^n}|\beta(s,X_{s-}^{n,t,x},e)-\beta(s,X_{s-}^{t,x},e)|^q\lambda(\mathrm de)\mathrm ds\\
            &\qquad\qquad\qquad\qquad\text{ } + C_q\mathbb{E}\left\{\left(\int^T_0\int_{E\setminus E^n}|\beta(s,X_{s-}^{t,x},e)|^2\lambda(\mathrm de)\mathrm ds\right)^{\frac{q}{2}}+ \int_0^T\int_{E\setminus E^n}|\beta(s,X_{s-}^{t,x},e)|^q\lambda(\mathrm de)\mathrm ds\right\}\\
            &\qquad\qquad\qquad\qquad \leq C_q\mathbb{E}\int_0^T|\Delta X^n_s|^q\mathrm ds + C_q \left\{\left(\int_{E\setminus E^n}\!(1\wedge |e|^2)\lambda(\mathrm de)\right)^{\frac{q}{2}}\!+\!\int_{E\setminus E^n}\!(1\wedge |e|^q)\lambda(\mathrm de)\right\}.
        \end{aligned}
    \end{equation*}
    The dominated convergence theorem gives $\int_{E\setminus E^n}(1\wedge |e|^2)\lambda(\mathrm de)\xrightarrow{n\to \infty}0$. Combining this with $\int_{E\setminus E^n}(1\wedge |e|^q)\lambda(\mathrm de)\leq \int_{E\setminus E^n}(1\wedge |e|^2)\lambda(\mathrm de)$ and Gronwall's inequality yields
    \begin{equation}
        \label{eq020}
        \lim_{n\to \infty} \mathbb{E} \left[\sup_{0\leq s\leq T}|\Delta X^n_s|^q\right] = 0,\qquad q\geq 2.
    \end{equation}

    We next prove convergence of the BSDE with Poisson jumps. Since $\mathbb{F}^n\subseteq \mathbb{F}$, \eqref{eq013} can be rewritten as the following BSDE with Poisson jumps under $\mathbb{F}$:
    \begin{equation}
        \label{eq021}
        \begin{aligned}
            &Y_{s}^{n,t,x}=g(X_T^{n,t,x})+\int_s^Tf\left(r,X_r^{n,t,x},Y_r^{n,t,x},Z_{r}^{n,t,x},\mathbbm{1}_{E^n}U_{r}^{n,t,x}\right)\mathrm dr\\
            &\qquad\qquad-\int_s^T Z_{r}^{n,t,x}\mathrm dW_r-\int_s^T\int_E\mathbbm{1}_{E^n}U_{r}^{n,t,x}\tilde\mu(\mathrm dr,\mathrm de),\qquad s\in [0,T].
        \end{aligned}
    \end{equation}
    To simplify notation, define $\Delta f^n_s:=f(s,X_s^{n,t,x},Y_s^{t,x},Z_s^{t,x},U_s^{t,x}) -f(s,X_s^{t,x},Y_s^{t,x},Z_s^{t,x},U_s^{t,x})$ and $\Delta g^n_T := g(X_T^{n,t,x})-g(X_T^{t,x})$. Since the constants in \Cref{re22}(iii) and \eqref{eq014} are independent of $n$, applying Proposition 2.2 of~\cite{BBP1997} to \eqref{eq004} and \eqref{eq021} gives
    \begin{equation}
        \label{eq024}
        \mathbb E\!\left[\sup_{0\le s\le T}|\Delta Y^n_s|^2+\int_0^T\!\Big(|\Delta Z^n_s|^2+\|\mathbbm{1}_{E^n}U_{s}^{n,t,x}-U_s^{t,x}\|_{L^2(E;\mathbb R^m)}^2\Big)\mathrm ds\right]\leq C\mathbb E\!\left[|\Delta g^n_T|^2+\int^T_0\!|\Delta f^n_s|^2\mathrm ds\right]. 
    \end{equation}
    Assumption~2(iv), \eqref{eq010}, \eqref{eq013_1}, and \eqref{eq020} imply
    \begin{equation}
        \label{eq025}
        \mathbb E\big[|\Delta g^n_T|^2\big]\leq C\left(\mathbb E\big[(1+|X_T^{n,t,x}|^p+|X_T^{t,x}|^p)^4\big]\right)^{\frac{1}{2}}\left(\mathbb E\big[|\Delta X^n_T|^4\big]\right)^{\frac{1}{2}}\xrightarrow{n\to \infty}0.
    \end{equation}
    Furthermore, by \eqref{eq007}, \eqref{eq010}, \eqref{eq013_1}, and \eqref{eq020},
    \begin{equation}
        \label{eq026}
        \begin{aligned}
            &\mathbb E\int^T_0|\Delta f^n_s|^2\mathrm ds \leq C\left(\mathbb E\left[\sup_{0\leq s\leq T}\big(1+|X_s^{n,t,x}|^p+|X_s^{t,x}|^p\big)^8\right]\right)^{\frac{1}{4}}\left(\mathbb E\left[\sup_{0\leq s\leq T}|\Delta X^n_s|^8\right]\right)^{\frac{1}{4}}\\
            &\qquad\qquad\qquad\qquad \times \left(\mathbb E\left[\left(T+\int_0^T\|U_s^{t,x}\|_{L^2(E;\mathbb R^m)}^2\mathrm ds\right)^2\right]\right)^{\frac{1}{2}}\xrightarrow{n\to \infty}0.
        \end{aligned}
    \end{equation}
    Combining \eqref{eq024}--\eqref{eq026}, we conclude that
    \begin{equation}
        \label{eq027}
        \lim_{n\to \infty} \mathbb E\left[\sup_{0\le s\le T}|\Delta Y^n_s|^2+\int_0^T\Big(|\Delta Z^n_s|^2+\|\mathbbm{1}_{E^n}U_{s}^{n,t,x}-U_s^{t,x}\|_{L^2(E;\mathbb R^m)}^2\Big)\mathrm ds\right] = 0.
    \end{equation}
{\bf \underline{Step 3}: Markovian representation of $U^{t,x}$.} By \eqref{eq006}, \eqref{eq016}, and \eqref{eq027}, for any $(t,x)\in [0,T]\times\mathbb R^{k}$,
    \begin{equation}
        \label{eq028}
        \lim_{n\to \infty} |u^n(t,x)-u(t,x)| = \lim_{n\to \infty} |\Delta Y^n_t| = 0.
    \end{equation}
    For any $(t,x)\in [0,T]\times\mathbb R^{k}$ and every sequence $\{x_n\}_{n\geq 1}\subseteq \mathbb R^{k}$ satisfying $\lim_{n\to \infty}x_n=x$, \eqref{eq018} and \eqref{eq028} imply
    \begin{equation}
        \label{eq029}
        \begin{aligned}
            &\lim_{n\to \infty}|u^n(t,x_n)-u(t,x)|\leq \lim_{n\to \infty}|u^n(t,x_n)-u^n(t,x)| + \lim_{n\to \infty}|u^n(t,x)-u(t,x)|\\
            &\qquad\qquad\qquad\qquad\qquad\text{ }\text{ }\leq \lim_{n\to \infty} C(1+|x_n|^p+|x|^p)|x_n-x| = 0.
        \end{aligned}
    \end{equation}
    Fix $(t,x)\in [0,T]\times\mathbb R^{k}$. It follows from \eqref{eq020} and \eqref{eq027} that the index sequence $\{n\}_{n\geq 1}$ admits an increasing subsequence $\{n_j\}_{j\geq 1}$ such that
    \begin{equation}
        \label{eq030}
        \left\{
            \begin{aligned}
                &\lim_{j\to \infty}\big|X_{s-}^{n_j,t,x}-X_{s-}^{t,x}\big|\leq \lim_{j\to \infty} \sup_{0\leq s\leq T}\big|\Delta X^{n_j}_s\big| = 0, \quad \forall s\in [0,T],\quad \mathbb{P}\text{-a.s.},\\
                &\lim_{j\to \infty}\big|\mathbbm{1}_{E^{n_j}}(e)U_{s}^{n_j,t,x}(e)-U_s^{t,x}(e)\big|=0,\quad \mathrm{d}\mathbb{P} \otimes \mathrm{d}\lambda \otimes \mathrm{d}t\text{-a.e. on } \Omega \times E \times [0,T].
            \end{aligned}
        \right.
    \end{equation}
    Hence, for any $(t,x)\in [0,T]\times\mathbb R^{k}$, Assumption~1(ii), \eqref{eq019}, \eqref{eq029}, and \eqref{eq030} yield
    \begin{equation}
        \label{eq032}
        \begin{aligned}
            & U^{t,x}_s(e) = \lim_{j\to \infty}\mathbbm{1}_{E^{n_j}}(e)U_{s}^{n_j,t,x}(e)= \lim_{j\to \infty}\mathbbm{1}_{E^{n_j}}(e)\left[u^{n_j}\big(s,X_{s-}^{n_j,t,x}+\beta(s,X_{s-}^{n_j,t,x},e)\big)-u^{n_j}\big(s,X_{s-}^{n_j,t,x}\big)\right]\\
            &\qquad\quad  = u\big(s,X_{s-}^{t,x}+\beta(s,X_{s-}^{t,x},e)\big)-u\big(s,X_{s-}^{t,x}\big),\quad \mathrm{d}\mathbb{P} \otimes \mathrm{d}\lambda \otimes \mathrm{d}t\text{-a.e. on } \Omega \times E \times [t,T].
        \end{aligned}
    \end{equation}
    This completes the proof.
\end{proof}

\begin{remark}
    \label{re33} 
    Some comments on the above proof are in order.
    \begin{enumerate}[(i)]
    \item In truncating the last argument of $h$ in \eqref{eq019_2}, the full per-jump contribution $\int_E\gamma_i(t,x,e)\Xi_i(0)\lambda(\mathrm de)$ is retained. If this term were replaced by $\int_E\gamma_i(t,x,e)\Xi_i(0)\lambda^n(\mathrm de)$ with $\Xi(0)\neq0$ and $\int_E(1\wedge|e|)\lambda(\mathrm de)=\infty$, Assumption~3(ii) would no longer apply directly to the resulting shifted generator. To establish the local Lipschitz continuity of $f^n$ with respect to $x$ using Assumption~2(i) and~(ii), the spatial variation of the truncated term would then need to be estimated separately. For any $x,x'\in\mathbb R^k$,
    \begin{equation*}
        \begin{aligned}
            &\quad \left|\int_E\gamma_i(t,x,e)\Xi_i(0)\lambda^n(\mathrm de)-\int_E\gamma_i(t,x',e)\Xi_i(0)\lambda^n(\mathrm de)\right|\\
            &\leq |\Xi_i(0)|\int_E|\gamma_i(t,x,e)-\gamma_i(t,x',e)|\lambda^n(\mathrm de)\\
            &\leq C_\gamma|\Xi(0)|\big(1+|x|^{p_\gamma}+|x'|^{p_\gamma}\big)|x-x'|\int_{E^n}(1\wedge|e|)\lambda(\mathrm de)\\
            &\leq nC_\gamma|\Xi(0)|\big(1+|x|^{p_\gamma}+|x'|^{p_\gamma}\big)|x-x'|\int_E(1\wedge|e|^2)\lambda(\mathrm de),
        \end{aligned}
    \end{equation*}
    where $C_\gamma$ and $p_\gamma$ are the constants in Assumption~2(ii). Consequently, this estimate does not provide a local Lipschitz constant $C$ independent of $n$ in the second inequality of \eqref{eq014}.
    \item In fact, to establish the Markovian representation of $U^{t,x}$ under an infinite L\'{e}vy measure, it suffices that the joint convergence \eqref{eq029} hold. Neither \eqref{eq018} nor \eqref{eq028} is required. This observation is crucial to the uniqueness proof of \Cref{th43}.
    \end{enumerate}
\end{remark}

\section{Well-Posedness of Viscosity Solutions to the IPDE}
\label{sec4}

We now present the main result of this paper, namely, the well-posedness of IPDE \eqref{eq001}. Recall that \eqref{eq001} takes the form
\begin{equation}\label{edpwp}
    \left\{
    \begin{aligned}
        &-\partial_tu_i(t,x)-\mathcal L u_i(t,x)- \mathcal Ku_i(t,x)\\
        &\qquad -h_i\left(t,x,u(t,x),\big(D_xu_i(t,x)\big)^\top \sigma(t,x),\mathcal B_i^\Xi u_i(t,x)\right)=0, \quad (t,x)\in[0,T]\times \mathbb{R}^k,\\
        & u_i(T,x)=g_i(x), \quad i = 1,\cdots,m,
    \end{aligned}
    \right.
\end{equation}
where
\begin{equation}
    \label{eq4_07}
    \left\{
        \begin{aligned}
            &\mathcal L u_i(t,x) = b(t,x)^\top D_xu_i(t,x)+\frac{1}{2}\operatorname{Tr}\left(\sigma(t,x)\sigma^\top(t,x)D^2_{xx}u_i(t,x)\right),\\
            &\mathcal K u_i(t,x) = \int_E\left[u_i(t,x+\beta(t,x,e))-u_i(t,x)-\beta(t,x,e)^\top D_xu_i(t,x)\right]\lambda(\mathrm{d}e),\\
            & \mathcal B_i^\Xi u_i(t,x) = \int_E\gamma_i(t,x,e)\Xi_i\left(u_i(t,x+\beta(t,x,e))-u_i(t,x)\right)\lambda(\mathrm{d}e).
        \end{aligned}
    \right.
\end{equation}
We use the following definition of a viscosity solution of equation \eqref{edpwp}.
\begin{definition}
    \label{de41}
    Let $u=(u_1,\cdots,u_m)^\top\in \mathcal U(\mathbb R^m)$ be a function. We call $u$ a viscosity sub-solution (resp. super-solution) of IPDE \eqref{eq001} if the following conditions hold:

    \begin{enumerate}[(i)]
        \item For any $x\in \mathbb R^k$, $u(T,x)\leq g(x)$ (resp. $u(T,x)\geq g(x)$);
        \item For any $i = 1,\cdots,m$ and $\varphi \in C^{1,2}([0,T]\times\mathbb R^k;\mathbb R)$, if $(t,x) \in [0,T)\times\mathbb R^k$ is a global maximum point (resp. minimum point) of $u_i-\varphi$ such that $u_i(t,x)-\varphi(t,x) = 0$, then
    \begin{equation*}
        -\partial_t\varphi(t,x)-\mathcal L \varphi(t,x)- \mathcal K\varphi(t,x)-h_i\left(t,x,u(t,x),\big(D_x\varphi(t,x)\big)^\top \sigma(t,x),\mathcal B_i^\Xi u_i(t,x)\right)\leq 0,
    \end{equation*}
    (resp.
    \begin{equation*}
        -\partial_t\varphi(t,x)-\mathcal L \varphi(t,x)- \mathcal K\varphi(t,x)-h_i\left(t,x,u(t,x),\big(D_x\varphi(t,x)\big)^\top \sigma(t,x),\mathcal B_i^\Xi u_i(t,x)\right)\geq 0,)
    \end{equation*}
    where $\mathcal L\varphi$ and $\mathcal K\varphi$ are defined by replacing $u_i$ with $\varphi$ in \eqref{eq4_07}. 
    \end{enumerate}
    If $u$ is both a viscosity sub-solution and a viscosity super-solution, it is called a viscosity solution of \eqref{eq001}.
\end{definition}

\begin{remark}
    \label{re42}
    Some comments on this definition are in order.
    \begin{enumerate}[(i)]
    \item Assumption~1(ii), Assumption~2(ii)--(iii), and Assumption~3(i) imply that $\mathcal B_i^\Xi \nu \in \Pi^c(\mathbb R)$ for any $\nu\in \mathcal U(\mathbb R)$. Thus, under the present assumptions, both IPDE \eqref{eq001} and \Cref{de41} are well defined.
    \item Unlike the classical definition of a viscosity solution in~\cite{BBP1997} (see \Cref{deA1}), \Cref{de41} assumes $u \in \mathcal U(\mathbb R^m)$ and retains $\mathcal B_i^\Xi u_i(t,x)$ in the viscosity inequalities. In Definition 3.1 and Lemma 3.3 of~\cite{BBP1997}, this term is replaced by $\mathcal B_i^\Xi \varphi(t,x)$, which is one of the main reasons for imposing monotonicity conditions on $\gamma$ and $h$ in~\cite{BBP1997}. If, however, $h_i$ does not depend on its last argument $\mathcal B_i^\Xi u_i(t,x)$, the viscosity inequalities in \Cref{de41} and \Cref{deA1} coincide.
    \end{enumerate}
\end{remark}

The main result is as follows.
\begin{theorem}
    \label{th43}
    Under Assumptions 1--3, the function $u$ in \eqref{eq006} is the unique viscosity solution of IPDE \eqref{eq001} in the sense of \Cref{de41}.
\end{theorem}

To prove \Cref{th43}, we first consider a special case of \eqref{eq001} in which the last argument $\mathcal B_i^\Xi u_i(t,x)$ of $h_i$ is replaced by a prescribed continuous field. More precisely, for $\vartheta = (\vartheta_1,\cdots,\vartheta_m)^\top \in \Pi^c(\mathbb R^m)$, consider the auxiliary IPDE
\begin{equation}
    \label{eq4_08}
    \left\{
    \begin{aligned}
        &-\partial_t\tilde u_i(t,x)-\mathcal L \tilde u_i(t,x)- \mathcal K\tilde u_i(t,x)\\
        &\qquad -h_i\left(t,x,\tilde u(t,x),\big(D_x\tilde u_i(t,x)\big)^\top \sigma(t,x),\vartheta_i(t,x)\right)=0, \quad (t,x)\in[0,T]\times \mathbb{R}^k,\\
        &\tilde u_i(T,x)=g_i(x), \quad i = 1,\cdots,m,
    \end{aligned}
    \right.
\end{equation}
where $\mathcal L\tilde u_i$ and $\mathcal K\tilde u_i$ are defined by replacing $u_i$ with $\tilde u_i$ in \eqref{eq4_07}. IPDE \eqref{eq4_08} falls within the special case for $h_i$ described in \Cref{re42}. The viscosity inequalities for \eqref{eq4_08} therefore coincide under \Cref{de41} and \Cref{deA1}, with the last argument of $h_i$ replaced by $\vartheta_i(t,x)$. The admissible solution class, however, is different: since $\vartheta \in \Pi^c(\mathbb R^m)$, viscosity solutions of \eqref{eq4_08} are sought in $\Pi^c(\mathbb R^m)$. As in \Cref{subse23}, for any $(t,x)\in[0,T]\times\mathbb R^k$, introduce the associated BSDE with Poisson jumps:
\begin{equation}
    \label{eq4_09}
    \tilde Y_{s}^{t,x}=g(X_T^{t,x})+\int_s^T\tilde f(r,X_r^{t,x},\tilde Y_r^{t,x},\tilde Z_{r}^{t,x})\mathrm dr-\int_s^T \tilde Z_{r}^{t,x}\mathrm dW_r-\int_s^T\int_E\tilde U_{r}^{t,x}\tilde\mu(\mathrm dr,\mathrm de),\quad s\in [0,T],
\end{equation}
where $\tilde f= (\tilde f_1,\cdots,\tilde f_m)^\top:[0,T]\times\mathbb R^k\times\mathbb R^m \times\mathbb R^{m\times d} \to\mathbb R^m$ is measurable, and its component functions $\tilde f_i:[0,T]\times\mathbb R^k\times\mathbb R^m \times\mathbb R^{1\times d}\to\mathbb R, i=1,\cdots,m$, are related to $\tilde f$ by
\begin{equation}
    \label{eq034}
    \begin{aligned}
        &\tilde f(t,x,y,z)=\big(\tilde f_1(t,x,y,z_1),\cdots,\tilde f_m(t,x,y,z_m)\big)^\top\\
        &=\left(h_1(t,x,y,z_1,\vartheta_1(t,x)),\cdots,h_m(t,x,y,z_m,\vartheta_m(t,x))\right)^\top,\text{ }\forall (t,x,y,z) \in [0,T]\times\mathbb R^k\times\mathbb R^m \times\mathbb R^{m\times d},
    \end{aligned}
\end{equation}
where $z_i$ is the $i$-th row of the matrix $z$. Under Assumption 1 and Assumption~2(i) and~(iv), Theorem 2.1 of~\cite{BBP1997} ensures that \eqref{eq4_09} admits a unique solution $(\tilde Y^{t,x},\tilde Z^{t,x},\tilde U^{t,x}) \in \mathcal S^2_{\mathbb F}([0,T];\mathbb R^m)\times \mathcal H^2_{\mathbb F}([0,T];\mathbb R^{m\times d}) \times \mathcal M^2_{\mathbb F}([0,T];\mathbb R^m)$. In view of \Cref{pro23}(ii) and \Cref{re24}(i), there exists a function $\tilde u=(\tilde u_1,\cdots,\tilde u_m)^\top\in \Pi^c(\mathbb R^m)$ such that, for any $(t,x)\in [0,T]\times\mathbb R^k$, $\tilde u(T,x) = g(x)$ and
\begin{equation}
    \label{eq4_10}
    \tilde Y_s^{t,x} = \tilde u(s,X_s^{t,x}),\quad \forall s\in [t,T],\quad \mathbb{P}\text{-a.s.} 
\end{equation}
Theorems 3.4 and 3.5 of~\cite{BBP1997} give the following well-posedness result for viscosity solutions of IPDE \eqref{eq4_08}.
\begin{lemma}
    \label{le44}
    Under Assumption 1 and Assumption~2(i) and~(iv), the function $\tilde u$ in \eqref{eq4_10} is the unique viscosity solution of IPDE \eqref{eq4_08}.
\end{lemma}

\begin{remark}
    \label{re45}
    Some comments on this lemma are in order.
    \begin{enumerate}[(i)]
    \item In IPDE \eqref{eq4_08}, the last argument of $h_i$ is the prescribed field $\vartheta_i$ and does not depend on $\gamma$ or the viscosity solution $\tilde u$. Hence, the monotonicity conditions (A.2iii) and (A.2iv) in~\cite{BBP1997} are not required.
    \item The condition $\vartheta \in \Pi^c(\mathbb R^m)$ and Assumption~2(i) ensure that assumption (A.2v) in~\cite{BBP1997} (related to uniqueness) holds. Indeed, for any $R>0$, define
    \begin{equation*}
        m_R(\rho) := C_h\left\{\big(1+2R^{p_h}\big)\rho + \sup_{\substack{t\in[0,T]\\|x|,|x'|\leq R,|x-x'|\leq \rho}}|\vartheta(t,x)-\vartheta(t,x')|\right\},\qquad \rho\geq 0,
    \end{equation*}
    where $C_h$ and $p_h$ are the constants in Assumption~2(i). The uniform continuity of $\vartheta$ on compact sets implies $\lim_{\rho\to 0}m_R(\rho) = 0$. Thus, for any $(t,x,x',y,z)\in[0,T]\times\mathbb R^{2k}\times\mathbb R^{m} \times\mathbb R^{1\times d}$ with $|x|,|x'|\leq R$,
    \begin{equation*}
        \begin{aligned}
            &|\tilde f_i(t,x,y,z)-\tilde f_i(t,x',y,z)|\leq |h_i\big(t,x,y,z,\vartheta_i(t,x)\big)-h_i\big(t,x',y,z,\vartheta_i(t,x)\big)|\\
            &\qquad\qquad\qquad\qquad\qquad\qquad\quad+|h_i\big(t,x',y,z,\vartheta_i(t,x)\big)-h_i\big(t,x',y,z,\vartheta_i(t,x')\big)|\\
            &\qquad\qquad\qquad\qquad\qquad\quad\text{ }\text{ }\leq m_R(|x-x'|),\qquad i = 1,\cdots,m.
        \end{aligned}
    \end{equation*}
    \item Although $b$, $\sigma$, and $\beta$ depend on time $t$, the continuity in $t$ required by Assumption 1 ensures that the errors due to this dependence vanish.
    \item In fact, uniqueness of the viscosity solution of IPDE \eqref{eq4_08} also holds in the class of continuous functions satisfying the following weaker growth condition:
    \begin{equation*}
        \lim_{|x|\to \infty}|\tilde u(t,x)|e^{-\tilde A \big[\ln(|x|)\big]^2} = 0,
    \end{equation*} 
    for some constant $\tilde A>0$, uniformly in $t\in[0,T]$. This class contains $\Pi^c(\mathbb R^m)$.
    \end{enumerate}
\end{remark}

We now turn to the proof of \Cref{th43}.
\begin{proof}[Proof of \Cref{th43}] {\bf \underline{Step 1}: Existence.} For any $(t,x)\in[0,T]\times\mathbb R^k$, define
    \begin{equation*}
        \vartheta(t,x) = (\vartheta_1(t,x),\cdots,\vartheta_m(t,x))^\top := (\mathcal B_1^\Xi u_1(t,x),\cdots,\mathcal B_m^\Xi u_m(t,x))^\top.
    \end{equation*}
    By \Cref{re42}(i), $\vartheta \in \Pi^c(\mathbb R^m)$. Since $X^{t,x}_{s-} = X^{t,x}_{s}, \mathrm{d}\mathbb{P} \otimes \mathrm{d}t$-a.e., \Cref{th32} gives
    \begin{equation*}
        \begin{aligned}
            &\int_E\gamma_i(s,X_s^{t,x},e)\Xi_i(U_{i,s}^{t,x}(e))\lambda(\mathrm de)=\int_E\gamma_i(s,X_{s-}^{t,x},e)\Xi_i(U_{i,s}^{t,x}(e))\lambda(\mathrm de)\\
            &\qquad\qquad\qquad\qquad\qquad\qquad\qquad=\int_E\gamma_i(s,X_{s-}^{t,x},e)\Xi_i\big(u_i(s,X_{s-}^{t,x}+\beta(s,X_{s-}^{t,x},e))-u_i(s,X_{s-}^{t,x})\big)\lambda(\mathrm de)\\
            &\qquad\qquad\qquad\qquad\qquad\qquad\qquad=\int_E\gamma_i(s,X_s^{t,x},e)\Xi_i\big(u_i(s,X_{s}^{t,x}+\beta(s,X_{s}^{t,x},e))-u_i(s,X_{s}^{t,x})\big)\lambda(\mathrm de)\\
            &\qquad\qquad\qquad\qquad\qquad\qquad\qquad=\vartheta_i(s,X_s^{t,x}), \qquad \mathrm{d}\mathbb{P} \otimes \mathrm{d}t\text{-a.e. on } \Omega \times [0,T],
        \end{aligned}
    \end{equation*}
    where $U_{i}^{t,x}$ is the $i$-th component of $U^{t,x}$. BSDE \eqref{eq004} can therefore be rewritten as
    \begin{equation*}
        Y_{s}^{t,x}=g(X_T^{t,x})+\int_s^T\tilde f(r,X_r^{t,x},Y_r^{t,x},Z_{r}^{t,x})\mathrm dr-\int_s^T Z_{r}^{t,x}\mathrm dW_r-\int_s^T\int_EU_{r}^{t,x}\tilde\mu(\mathrm dr,\mathrm de),\quad s\in [0,T],
    \end{equation*}
    where the generator $\tilde f$ is defined in \eqref{eq034}. \Cref{le44} implies that, for any $i = 1,\cdots,m$ and $\varphi \in C^{1,2}([0,T]\times\mathbb R^k;\mathbb R)$, if $(t,x) \in [0,T)\times\mathbb R^k$ is a global maximum point of $u_i-\varphi$ with $u_i(t,x)-\varphi(t,x) = 0$, then
    \begin{equation*}
        \begin{aligned}
            &\quad -\partial_t\varphi(t,x)-\mathcal L \varphi(t,x)- \mathcal K\varphi(t,x)-h_i\left(t,x,u(t,x),\big(D_x\varphi(t,x)\big)^\top \sigma(t,x),\mathcal B_i^\Xi u_i(t,x)\right)\\
            &= -\partial_t\varphi(t,x)-\mathcal L \varphi(t,x)- \mathcal K\varphi(t,x)-h_i\left(t,x,u(t,x),\big(D_x\varphi(t,x)\big)^\top \sigma(t,x),\vartheta_i(t,x)\right)\leq 0.
        \end{aligned}
    \end{equation*}
    Thus, $u$ is a viscosity sub-solution of \eqref{eq001}. A similar argument shows that $u$ is a viscosity super-solution of \eqref{eq001}. Consequently, $u$ is a viscosity solution of IPDE \eqref{eq001} in the sense of \Cref{de41}.

\noindent {\bf \underline{Step 2}: Uniqueness.} Let $\bar u = (\bar u_1,\cdots,\bar u_m)^\top \in \mathcal{U}(\mathbb{R}^m)$ be another viscosity solution of IPDE \eqref{eq001} in the sense of \Cref{de41}. Set $\bar \vartheta = (\bar \vartheta_1,\cdots,\bar \vartheta_m)^\top := (\mathcal B_1^\Xi \bar u_1,\cdots,\mathcal B_m^\Xi \bar u_m)^\top$. By \Cref{re42}(i), $\bar\vartheta\in\Pi^c(\mathbb R^m)$. Consider now the following IPDE:
    \begin{equation}
        \label{eq036}
        \left\{
            \begin{aligned}
                &-\partial_t\tilde u_i(t,x)-\mathcal L \tilde u_i(t,x)- \mathcal K\tilde u_i(t,x)\\
                &\qquad -h_i\left(t,x,\tilde u(t,x),\big(D_x\tilde u_i(t,x)\big)^\top \sigma(t,x),\bar \vartheta_i(t,x)\right)=0, \quad (t,x)\in[0,T]\times \mathbb{R}^k,\\
                &\tilde u_i(T,x)=g_i(x), \quad i = 1,\cdots,m.
            \end{aligned}
        \right.
    \end{equation}
    The BSDE with Poisson jumps associated with \eqref{eq036} is
    \begin{equation}
        \label{eq037}
        \bar Y_{s}^{t,x}=g(X_T^{t,x})+\int_s^T\bar  f(r,X_r^{t,x},\bar  Y_r^{t,x},\bar  Z_{r}^{t,x})\mathrm dr-\int_s^T \bar  Z_{r}^{t,x}\mathrm dW_r-\int_s^T\int_E\bar  U_{r}^{t,x}\tilde\mu(\mathrm dr,\mathrm de),\quad s\in [0,T],
    \end{equation}
    where, for any $(t,x,y)\in[0,T]\times\mathbb R^{k}\times\mathbb R^{m}$ and $z = (z_1,\cdots,z_m)^\top \in \mathbb R^{m\times d}$,
    \begin{equation*}
        \bar f(t,x,y,z) \!=\! \big(\bar f_1(t,x,y,z_1),\cdots\!,\bar f_m(t,x,y,z_m)\big)^\top\!=\!\left(h_1(t,x,y,z_1,\bar \vartheta_1(t,x)),\cdots\!,h_m(t,x,y,z_m,\bar \vartheta_m(t,x))\right)^\top.
    \end{equation*}
    Theorem 2.1 of~\cite{BBP1997} guarantees that \eqref{eq037} admits a unique solution $(\bar Y^{t,x},\bar Z^{t,x},\bar U^{t,x}) \in \mathcal S^2_{\mathbb F}([0,T];\mathbb R^m)\times \mathcal H^2_{\mathbb F}([0,T];\mathbb R^{m\times d}) \times \mathcal M^2_{\mathbb F}([0,T];\mathbb R^m)$. Furthermore, \Cref{pro23}(ii) and \Cref{re24}(i) provide a deterministic function $\tilde u\in \Pi^c(\mathbb R^m)$ such that, for any $(t,x)\in [0,T]\times\mathbb R^k$, $\tilde u(T,x) = g(x)$ and
    \begin{equation}
        \label{eq038}
        \bar Y_s^{t,x} = \tilde u(s,X_s^{t,x}),\quad \forall s\in [t,T],\quad \mathbb{P}\text{-a.s.} 
    \end{equation}
    By \Cref{le44}, $\tilde u$ is the unique viscosity solution of \eqref{eq036}. Moreover, $\bar u$ is also a viscosity solution of the same IPDE. Indeed, $\bar u\in\mathcal U(\mathbb R^m)\subseteq\Pi^c(\mathbb R^m)$, and since $\bar\vartheta_i=\mathcal B_i^\Xi\bar u_i$, conditions (i) and (ii) of \Cref{de41} for $\bar u$ coincide with those for \eqref{eq036}. Hence, $\bar u=\tilde u$, and \eqref{eq038} gives
    \begin{equation}
        \label{eq039}
        \bar Y_s^{t,x} = \tilde u(s,X_s^{t,x}) = \bar u(s,X_s^{t,x}),\quad \forall s\in [t,T],\quad \mathbb{P}\text{-a.s.}
    \end{equation}
    We next establish a representation of $\bar U^{t,x}$ in terms of $\bar u$ and $X^{t,x}$. Since $\bar \vartheta \in \Pi^c(\mathbb R^m)$, the function $\bar f(\cdot,\cdot,y,z)$ need not belong to $\mathcal{U}(\mathbb R^m)$, so $\bar f$ may fail to satisfy Assumption~2(i), and \Cref{th32} cannot be applied directly. Nevertheless, we follow the same proof strategy, first truncating the L\'{e}vy measure to a finite measure and then passing to the limit. In view of \Cref{re33}(ii), the difference is that we establish the joint convergence \eqref{eq029} directly.

    Recall the truncated Poisson random measure in \eqref{eq014_1}. For any $(t,x)\in [0,T]\times\mathbb R^{k}$ and any sequence $\{x_n\}_{n\geq 1}\subseteq \mathbb R^{k}$ satisfying $\lim_{n\to \infty}x_n=x$, consider the SDE driven by this measure:
    \begin{equation}
        \label{eq042}
        \left\{
            \begin{aligned}
                &X_s^{n,t,x_n}\!=\!x_n\!+\!\int_t^s \!b(r,X_r^{n,t,x_n})\mathrm dr \!+\!\int_t^s\!\sigma(r,X_r^{n,t,x_n})\mathrm dW_r\!+\!\int_t^s\!\int_E\!\beta(r,X_{r-}^{n,t,x_n},e)\tilde\mu^n(\mathrm dr,\mathrm de),\text{ }\text{ }s\in(t,T],\\
                &X_s^{n,t,x_n}=x_n,\qquad s\in[0,t].
            \end{aligned}
        \right.
    \end{equation}
    For any $n\geq 1$, Theorem 2.1 of~\cite{FK1985} gives a unique adapted solution $X^{n,t,x_n}$ to \eqref{eq042}. As in \eqref{eq013_1},
    \begin{equation}
        \label{eq043}
        \sup_{n\ge1}\mathbb E\left[\sup_{0\le s\le T}|X_s^{n,t,x_n}|^q\right]\le C_q\left(1+\sup_{n\geq 1}|x_n|^q\right)< \infty,\quad q\geq 2,
    \end{equation}
    where $C_q>0$ is a constant. For \eqref{eq002} and \eqref{eq042}, the estimates \eqref{eq013_1} and \eqref{eq020} yield
    \begin{equation}
        \label{eq044_1}
        \begin{aligned}
            & \quad \text{ }\lim_{n\to \infty} \mathbb{E}\left[\sup_{0\leq s\leq T}|X_s^{n,t,x_n}-X_s^{t,x}|^q\right]\\
            &\leq \lim_{n\to \infty} 2^{q-1}\mathbb{E}\left[\sup_{0\leq s\leq T}|X_s^{n,t,x_n}-X_s^{n,t,x}|^q\right]+\lim_{n\to \infty} 2^{q-1}\mathbb{E}\left[\sup_{0\leq s\leq T}|X_s^{n,t,x}-X_s^{t,x}|^q\right] = 0,\text{ }\text{ }q\geq 2.
        \end{aligned}
    \end{equation}
    Since $\bar \vartheta \in \Pi^c(\mathbb R^m)$, \eqref{eq044_1} and Theorem 6 in Chapter 2 of~\cite{GS1969} imply that $\bar \vartheta(\cdot,X_\cdot^{n,t,x_n})$ converges to $\bar \vartheta(\cdot,X_\cdot^{t,x})$ in measure with respect to the product measure $\mathrm{d}\mathbb{P} \otimes \mathrm{d}t$. Moreover, the family $\{|\bar \vartheta(\cdot,X_\cdot^{n,t,x_n})-\bar \vartheta(\cdot,X_\cdot^{t,x})|^2\}_{n\geq 1}$ is readily seen to be uniformly integrable on $\Omega\times[0,T]$. The Lebesgue--Vitali theorem (see Theorem 4.5.4 of~\cite{BR2007}) then gives
    \begin{equation}
        \label{eq045_1}
        \lim_{n\to \infty} \mathbb{E}\int^T_0\big|\bar \vartheta(s,X_s^{n,t,x_n})-\bar \vartheta(s,X_s^{t,x})\big|^2\mathrm ds = 0.
    \end{equation}
    Next, consider the following BSDE with the truncated Poisson random measure, corresponding to \eqref{eq037}:
    \begin{equation}
        \label{eq044}
        \begin{aligned}
            &\bar Y_{s}^{n,t,x_n}=g(X_T^{n,t,x_n})+\int_s^T\bar f\left(r,X_r^{n,t,x_n},\bar Y_r^{n,t,x_n},\bar Z_{r}^{n,t,x_n}\right)\mathrm dr-\int_s^T \bar Z_{r}^{n,t,x_n}\mathrm dW_r-\!\int_s^T\!\!\int_E\bar U_{r}^{n,t,x_n}\tilde\mu^n(\mathrm dr,\mathrm de)\\
            &\qquad\quad=  g(X_T^{n,t,x_n})+\int_s^T\bar f\left(r,X_r^{n,t,x_n},\bar Y_r^{n,t,x_n},\bar Z_{r}^{n,t,x_n}\right)\mathrm dr-\int_s^T \bar Z_{r}^{n,t,x_n}\mathrm dW_r\\
            &\qquad\qquad\quad-\int_s^T\int_E\mathbbm{1}_{E^n}\bar U_{r}^{n,t,x_n}\tilde\mu(\mathrm dr,\mathrm de),\qquad s\in [0,T].
        \end{aligned}
    \end{equation}
    For any $n\geq 1$, Theorem 2.1 of~\cite{BBP1997} yields a unique adapted solution $(\bar Y^{n,t,x_n},\bar Z^{n,t,x_n},\bar U^{n,t,x_n})$ of \eqref{eq044}. Combining \Cref{pro23}(ii) with \Cref{re24}(i) gives a deterministic function $\bar u^n\in \Pi^c(\mathbb R^m)$ such that, for any $(t,x)\in [0,T]\times\mathbb R^k$, $\bar u^n(T,x) = g(x)$ and
    \begin{equation}
        \label{eq045}
        \bar Y_s^{n,t,x_n} = \bar u^n(s,X_s^{n,t,x_n}),\quad \forall s\in [t,T],\quad \mathbb{P}\text{-a.s.} 
    \end{equation}
    In addition, \Cref{pro25}, \Cref{re26}, and \eqref{eq045} imply
    \begin{equation}
        \label{eq046}
       \mathbbm{1}_{E^{n}}(e) \bar U^{n,t,x_n}_s(e) =\bar  u^n\big(s,X_{s-}^{n,t,x_n}+\beta(s,X_{s-}^{n,t,x_n},e)\big)-\bar u^n\big(s,X_{s-}^{n,t,x_n}\big),\text{ }\text{ } \mathrm{d}\mathbb{P} \otimes \mathrm{d}\lambda \otimes \mathrm{d}t\text{-a.e. on } \Omega \times E \times [t,T].
    \end{equation}
    Applying Proposition 2.2 of~\cite{BBP1997} to \eqref{eq037} and \eqref{eq044} yields
    \begin{equation*}
        \begin{aligned}
            &\quad \text{ }\mathbb E\left[\sup_{0\le s\le T}|\bar Y_{s}^{n,t,x_n}-\bar Y_{s}^{t,x}|^2+\int_0^T\Big(|\bar Z_{s}^{n,t,x_n}-\bar Z_{s}^{t,x}|^2+\|\mathbbm{1}_{E^n}\bar U_{s}^{n,t,x_n}-\bar U_s^{t,x}\|_{L^2(E;\mathbb R^m)}^2\Big)\mathrm ds\right]\\
            &\leq C\mathbb E\left[|g(X_{T}^{n,t,x_n})-g(X_{T}^{t,x})|^2+\int^T_0|\bar f(s,X_s^{n,t,x_n},\bar Y_s^{t,x},\bar Z_s^{t,x}) -\bar f(s,X_s^{t,x},\bar Y_s^{t,x},\bar Z_s^{t,x})|^2\mathrm ds\right]\\
            &\leq CC^2_g\mathbb E\left[\big(1+|X_{T}^{n,t,x_n}|^{p_g}+|X_{T}^{t,x}|^{p_g}\big)^2\big|X_{T}^{n,t,x_n}-X_{T}^{t,x}\big|^2\right]\\
            &\text{ }\text{ }+2CC^2_h\mathbb E\int^T_0\left[\big(1+|X_{s}^{n,t,x_n}|^{p_h}+|X_{s}^{t,x}|^{p_h}\big)^2\big|X_{s}^{n,t,x_n}-X_{s}^{t,x}\big|^2+\big|\bar \vartheta(s,X_s^{n,t,x_n})-\bar \vartheta(s,X_s^{t,x})\big|^2\right]\mathrm ds,
        \end{aligned}
    \end{equation*}
    where the constant $C>0$ is independent of $n$, and $C_h,C_g,p_h,p_g$ are as in Assumption~2(i) and~(iv). Combining this estimate with \eqref{eq010} and \eqref{eq043}--\eqref{eq045_1}, we obtain
    \begin{equation}
        \label{eq049}
        \lim_{n\to\infty} \mathbb E\left[\sup_{0\le s\le T}|\bar Y_{s}^{n,t,x_n}-\bar Y_{s}^{t,x}|^2+\int_0^T\Big(|\bar Z_{s}^{n,t,x_n}-\bar Z_{s}^{t,x}|^2+\|\mathbbm{1}_{E^n}\bar U_{s}^{n,t,x_n}-\bar U_s^{t,x}\|_{L^2(E;\mathbb R^m)}^2\Big)\mathrm ds\right]  = 0.
    \end{equation}
    It follows from \eqref{eq039}, \eqref{eq045}, and \eqref{eq049} that the convergence property \eqref{eq029} holds: for any $(t,x)\in [0,T]\times\mathbb R^{k}$ and every sequence $\{x_n\}_{n\geq 1}\subseteq \mathbb R^{k}$ satisfying $\lim_{n\to \infty}x_n=x$,
    \begin{equation}
        \label{eq050}
        \lim_{n\to \infty} |\bar u^n(t,x_n)-\bar u(t,x)| = \lim_{n\to \infty} |\bar Y_{t}^{n,t,x_n}-\bar Y_{t}^{t,x}| = 0.
    \end{equation}
    Using \eqref{eq046}--\eqref{eq050} and arguing as in the proofs of \eqref{eq030} and \eqref{eq032}, we conclude that
    \begin{equation}
        \label{eq051}
        \bar U^{t,x}_s(e) = \bar u\big(s,X_{s-}^{t,x}+\beta(s,X_{s-}^{t,x},e)\big)-\bar u\big(s,X_{s-}^{t,x}\big),\quad \mathrm{d}\mathbb{P} \otimes \mathrm{d}\lambda \otimes \mathrm{d}t\text{-a.e. on } \Omega \times E \times [t,T].
    \end{equation}

    To complete the uniqueness proof, use \eqref{eq051} to rewrite \eqref{eq037} as
    \begin{equation*}
        \bar Y_{s}^{t,x}=g(X_T^{t,x})+\int_s^T f(r,X_r^{t,x},\bar  Y_r^{t,x},\bar  Z_{r}^{t,x},\bar  U_{r}^{t,x})\mathrm dr-\int_s^T \bar  Z_{r}^{t,x}\mathrm dW_r-\int_s^T\int_E\bar  U_{r}^{t,x}\tilde\mu(\mathrm dr,\mathrm de),\quad s\in [0,T],
    \end{equation*}
    where $f$ is defined in \eqref{eq004_1}. Thus, \eqref{eq004} and \eqref{eq037} are the same BSDE with Poisson jumps. For any $(t,x)\in [0,T]\times\mathbb R^{k}$, uniqueness of the solution to \eqref{eq004}, together with \eqref{eq006} and \eqref{eq039}, implies
    \begin{equation*}
        u(t,x) = Y^{t,x}_t = \bar Y^{t,x}_t = \bar u(t,x).
    \end{equation*}
    This completes the proof.
\end{proof}

\section{Conclusion}
\label{sec5}
This paper has studied viscosity solutions of second-order IPDE systems with nonlinear jump terms. In contrast to existing results that rely on monotonicity conditions~\cite{BBP1997} or involve only linear integrals of jump increments~\cite{H2016}, our framework requires $\Xi$ to be only Lipschitz continuous, allows its components to be nonmonotone with $\Xi(0)\neq0$, and imposes no monotonicity conditions on the system. We construct the solution function $u$ of the IPDE through BSDEs with Poisson jumps. Since the L\'{e}vy measure may be infinite, it is first truncated near the origin to obtain finite measures. Combining the local boundedness of the spatial difference quotients of $u$ with convergence of the truncated equations yields the Markovian representation of the jump component $U^{t,x}$. We then introduce an auxiliary IPDE driven by a prescribed continuous field $\vartheta$ and use the joint convergence of the functions associated with the truncated equations with respect to the truncation parameter and the initial state to prove that $u$ is the unique viscosity solution of the IPDE \eqref{eq001}. However, when $\Xi$ no longer satisfies the Lipschitz condition in Assumption~2(iii), the generator of the associated BSDE with Poisson jumps may fail to be Lipschitz continuous with respect to the jump component, and both the well-posedness of this BSDE and the spatial regularity of $u$ need to be established anew. Therefore, extending the results to the case where $\Xi$ is not Lipschitz continuous is a problem worthy of consideration, which we leave for future research.

\begin{appendix}
\section{Classical Definition of a Viscosity Solution}
    \label[appendix]{appA}
    For the comparison in \Cref{re42} between \Cref{de41} and the classical notion of a viscosity solution, we recall below the definition given in Definition 3.1 and Lemma 3.3 of~\cite{BBP1997}.

    \begin{definition}
        \label{deA1}
        Let $u=(u_1,\cdots,u_m)^\top:[0,T]\times\mathbb R^k\to\mathbb R^m$ be a function that is continuous on $[0,T]\times\mathbb R^k$. We call $u$ a viscosity sub-solution (resp. super-solution) of IPDE \eqref{eq001} if the following conditions hold:

        \begin{enumerate}[(i)]
            \item For any $x\in \mathbb R^k$, $u(T,x)\leq g(x)$ (resp. $u(T,x)\geq g(x)$);
            \item For any $i = 1,\cdots,m$ and $\varphi \in C^{1,2}([0,T]\times\mathbb R^k;\mathbb R)$, if $(t,x) \in [0,T)\times\mathbb R^k$ is a global maximum point (resp. minimum point) of $u_i-\varphi$ such that $u_i(t,x)-\varphi(t,x) = 0$, then
        \begin{equation*}
            -\partial_t\varphi(t,x)-\mathcal L \varphi(t,x)- \mathcal K\varphi(t,x)-h_i\left(t,x,u(t,x),\big(D_x\varphi(t,x)\big)^\top \sigma(t,x),\mathcal B_i^\Xi \varphi(t,x)\right)\leq 0,
        \end{equation*}
        (resp.
        \begin{equation*}
            -\partial_t\varphi(t,x)-\mathcal L \varphi(t,x)- \mathcal K\varphi(t,x)-h_i\left(t,x,u(t,x),\big(D_x\varphi(t,x)\big)^\top \sigma(t,x),\mathcal B_i^\Xi \varphi(t,x)\right)\geq 0,)
        \end{equation*}
        where $\mathcal L\varphi$, $\mathcal K\varphi$, and $\mathcal B_i^\Xi\varphi$ are defined by replacing $u_i$ with $\varphi$ in \eqref{eq4_07}. 
        \end{enumerate}
        If $u$ is both a viscosity sub-solution and a viscosity super-solution, it is called a viscosity solution of \eqref{eq001}.
    \end{definition}
\end{appendix}


\begin{thebibliography}{99}

\bibitem{BBP1997} Barles G, Buckdahn R, Pardoux E. Backward stochastic differential equations and integral-partial differential equations. Stochastics Stochastics Rep, 1997, 60: 57--83

\bibitem{B1997} Basu S. The conservatism principle and the asymmetric timeliness of earnings. J Account Econ, 1997, 24: 3--37

\bibitem{B2006} Becherer D. Bounded solutions to backward SDE's with jumps for utility optimization and indifference hedging. Ann Appl Probab, 2006, 16: 2027--2054

\bibitem{BHK2023} Bensalem S, Hern\'{a}ndez-Santib\'{a}\~{n}ez N, Kazi-Tani N. A continuous-time model of self-protection. Finance Stoch, 2023, 27: 503--537


\bibitem{BLY2022} Bo L, Liao H, Yu X. Risk-sensitive credit portfolio optimization under partial information and contagion risk. Ann Appl Probab, 2022, 32: 2355--2399

\bibitem{BR2007} Bogachev V I. Measure Theory. Berlin: Springer, 2007

\bibitem{BCCS2024} Brachetta M, Callegaro G, Ceci C, et al. Optimal reinsurance via BSDEs in a partially observable model with jump clusters. Finance Stoch, 2024, 28: 453--495

\bibitem{BHL2011} Buckdahn R, Hu Y, Li J. Stochastic representation for solutions of Isaacs' type integral--partial differential equations. Stochastic Process Appl, 2011, 121: 2715--2750

\bibitem{DQS2015} Dumitrescu R, Quenez M C, Sulem A. Optimal stopping for dynamic risk measures with jumps and obstacle problems. J Optim Theory Appl, 2015, 167: 219--242

\bibitem{DRZ2021} Dumitrescu R, Reisinger C, Zhang Y. Approximation schemes for mixed optimal stopping and control problems with nonlinear expectations and jumps. Appl Math Optim, 2021, 83: 1387--1429

\bibitem{EFO2017} Eddahbi M, Fakhouri I, Ouknine Y. $\mathbb{L}^p(p\geq 2)$-solutions of generalized BSDEs with jumps and monotone generator in a general filtration. Mod Stoch Theory Appl, 2017, 4: 25--63


\bibitem{FK1985} Fujiwara T, Kunita H. Stochastic differential equations of jump type and L\'{e}vy processes in diffeomorphisms group. J Math Kyoto Univ, 1985, 25: 71--106

\bibitem{FS2018} Funke B, Schmisser \'{E}. Adaptive nonparametric drift estimation of an integrated jump diffusion process. ESAIM Probab Stat, 2018, 22: 236--260
\bibitem{GS1969} Gikhman I I, Skorokhod A V. Introduction to the Theory of Random Processes. Philadelphia: W. B. Saunders, 1969

\bibitem{H2016} Hamad\`{e}ne S. Viscosity solutions of second order integral--partial differential equations without monotonicity condition: a new result. Nonlinear Anal, 2016, 147: 213--235

\bibitem{HMN2023} Hamad\`{e}ne S, Mnif M, Neffati S. Viscosity solution of system of integro-partial differential equations with interconnected obstacles of non-local type without monotonicity conditions. J Dynam Differential Equations, 2023, 35: 1151--1173

\bibitem{HMN2026} Hamad\`{e}ne S, Mnif M, Neffati S. Viscosity solution of systems of integral-partial differential equations with interconnected obstacles without monotonicity conditions and infinite L\'{e}vy measure. arXiv:2609.16812, 2026

\bibitem{HM2016} Hamad\`{e}ne S, Morlais M A. Viscosity solutions for second order integro-differential equations without monotonicity condition: the probabilistic approach. Stochastics, 2016, 88: 632--649

\bibitem{JP2011} Jacod J, Protter P. Discretization of Processes. Berlin: Springer, 2012

\bibitem{KPZ2015} Kazi-Tani N, Possama\"{\i} D, Zhou C. Quadratic BSDEs with jumps: a fixed-point approach. Electron J Probab, 2015, 20: 66, 28pp

\bibitem{LW2026} Liu L, Wei Q. Reflected stochastic recursive control problems with jumps: dynamic programming and stochastic verification theorems. ESAIM Control Optim Calc Var, 2026, 32: 36, 54pp

\bibitem{LLLW2025} Luo S, Li W, Li X, et al. Stochastic representation for solutions of a system of coupled HJB-Isaacs equations with integral--differential operators. Stochastic Process Appl, 2025, 179: 104502, 25pp

\bibitem{LLW2025} Luo S, Li X, Wei Q. Infinite time horizon stochastic recursive control problems with jumps: dynamic programming and stochastic verification theorems. SIAM J Control Optim, 2025, 63: 796--821

\bibitem{M1974} Merton R C. On the pricing of corporate debt: the risk structure of interest rates. J Finance, 1974, 29: 449--470

\bibitem{M2010} Morlais M A. A new existence result for quadratic BSDEs with jumps with application to the utility maximization problem. Stochastic Process Appl, 2010, 120: 1966--1995

\bibitem{OPK2023} Orthaus S, Pelger C, Kuhner C. The eternal debate over conservatism and prudence: a historical perspective on the conceptualization of asymmetry in financial accounting theory. Contemp Account Res, 2023, 40: 41--88

\bibitem{QS2013} Quenez M C, Sulem A. BSDEs with jumps, optimization and applications to dynamic risk measures. Stochastic Process Appl, 2013, 123: 3328--3357

\bibitem{SST2024} Santacroce M, Siri P, Trivellato B. Forward backward SDEs systems for utility maximization in jump diffusion models. Appl Math Optim, 2024, 89: 65, 22pp

\bibitem{W2003} Watts R L. Conservatism in accounting part I: explanations and implications. Account Horiz, 2003, 17: 207--221


\end{thebibliography}
\end{document}